\documentclass[11pt, twoside, fleqno]{article} 

\usepackage{stmaryrd}
\usepackage{amsmath}
\usepackage{a4}
\usepackage{euscript}
\usepackage{amsfonts}
\usepackage{amssymb}
\usepackage{mathrsfs}
\usepackage{xypic}
\usepackage{enumerate}
\usepackage{array}
\usepackage{theorem}
\xyoption{all}

\DeclareMathAlphabet{\mathpzc}{OT1}{pzc}{m}{it}

\newtheorem{theorem}{Theorem}[section]
\newtheorem{proposition}[theorem]{Proposition}
\newtheorem{corollary}[theorem]{Corollary}
\newtheorem{lemma}[theorem]{Lemma}
\theorembodyfont{\upshape}
\newtheorem{definition}[theorem]{Definition}
\newtheorem{remark}[theorem]{Remark}

\numberwithin{equation}{section}

\newcommand{\qed}{\hfill $\Box$ \medskip}

\newcommand{\NN}{{\mathbb N}}
\newcommand{\ZZ}{{\mathbb Z}}
\newcommand{\QQ}{{\mathbb Q}}

\newcommand{\PP}{{\mathbb P}}

\newcommand{\FF}{{\mathbb F}}

\renewcommand{\O}{{\cal O}}
\renewcommand{\L}{{\cal L}}

\newcommand{\Y}{{\cal Y}}

\newcommand{\GG}{\ensuremath{\mathbb G}}

\newcommand{\Spec}{{\operatorname{Spec}\kern 1pt}}
\newcommand{\Spf}{{\operatorname{Spf}\kern 1pt}}

\newcommand{\ord}{{\mathrm{ord}}}

\newcommand{\Ext}{{\mathrm{Ext}}}

\newcommand{\Hom}{{\mathrm{Hom}}}
\renewcommand{\hom}{{\mathpzc{Hom}}}
\renewcommand{\H}{{\mathrm{H}}}

\renewcommand{\Im}{\operatorname{Im}\kern 1pt}
\newcommand{\coker}{\operatorname{coker}\kern 1pt}

\renewcommand{\epsilon}{\varepsilon}

\newcommand{\Frac}{\mathrm{Frac}}

\newcommand{\proof}{\noindent{\it Proof : }}

\newcommand{\Som}{{\mathbf{Som}}}
\newcommand{\Ar}{{\mathbf{Ar}\kern 0.5pt}}
\newcommand{\supp}{{\mathrm{supp}}}

\newcommand{\limproj}{\mathop{\mathop\mathrm {lim}\limits_{\longleftarrow}}}

\newcommand{\X}{{\mathcal X}}

\newcommand{\U}{{\mathcal U}}

\newcommand{\p}{{\mathfrak p}}

\newcommand{\A}{{\mathfrak A}}

\newcommand{\m}{{\mathfrak m}}

\newcommand{\F}{{\mathcal F}}

\newcommand{\dif}{{\mathfrak d}}

\begin{document}

\newcommand{\Sepa}{\mathpzc{Sepa}}
\newcommand{\Defsep}{\mathpzc{DefSep}}
\newcommand{\af}{\mathfrak{a}}
\newcommand{\ca}{\mathrm{C}}
\newcommand{\Mod}{\mathrm{Mod}}
\newcommand{\Pp}{\mathscr{P}}

\title{On a compactification of a Hurwitz space in the wild case}

\author{Sylvain Maugeais}

\date{}

\maketitle

\begin{abstract}
Let $g$ and $g'$ be two integers and $p$ a prime number. Denote
by $\mathscr H_{g, g', p}^c$ the moduli space of morphisms of degree $p$
between smooth curves of genus $g$ and $g'$ and
with constant ramification. The purpose of this 
article is to construct and describe a modular compactification of 
this space. 
\end{abstract} 

\section{Introduction}

The aim of this work is to construct a compactification of one
Hurwitz space in the wild case. By Hurwitz spaces, we mean 
spaces classifying finite morphisms between curves with some 
condition (for example morphisms of a fixed degree, morphisms to 
$\PP^1$, ...). In the tame case, we already know a good modular 
compactification of those spaces (see for example the work
of Harris and Mumford \cite{Harris_Mumford}, Mochizuki 
\cite{Mochizuki} or Wewers \cite{TheseWewers}).
On the other hand, we are far from having a good compactification
in the wild case. Let us illustrate the problems that can arise
on an example. Let $g, g', d \in \NN$ with $g \ge 2$ and consider the 
moduli space ${\mathscr H}_{g, g', d}$ of geometrically separable 
morphisms of degree $d$ between proper smooth curves of genus $g$ 
and $g'$. Then it is relatively easy to prove that 
${\mathscr H}_{g, g', d}$ is representable by an algebraic stack that
is formally smooth over $\ZZ$.
Using stable curves, it is possible to embed this space in
a proper space over $\ZZ$. Indeed, consider the moduli space
$\bar {\mathscr H}_{g, g', d}$ of ``stable'' morphisms of degree $d$ 
between semi-stable curves of genus $g$ and $g'$ (the ``stable'' 
condition is just a minimality condition that ensures the space 
$\bar {\mathscr H}_{g, g', d}$ to be separated, cf. 
\cite{LiuStableReduction}). In general the space 
$\bar {\mathscr H}_{g, g', d}$ contains morphisms between singular curves
 and inseparable
morphisms (through a ``stable reduction''-like theorem, we can
see that this phenomenon is necessary in order to have a
proper space).
Then it is possible to prove that this space is representable by
a proper algebraic stack over $\ZZ$. The problems are that
in general, ${\mathscr H}_{g, g', d}$ fails to be dense in 
$\bar {\mathscr H}_{g, g', d}$ and that $\bar {\mathscr H}_{g, g', d}$ may
have nasty singularities. Moreover, it does not seem easy
to give a thorough description of the geometric points of the
closure of ${\mathscr H}_{g, g', d}$ (in particular, this closure
have probably no ``simple'' interpretation).

Our aim is to construct a compactification with mild singularities
of a variation of ${\mathscr H}_{g, g', p}$ for a prime $p$. The 
restriction on $p$ is not just a convenience matter. There are 
many instances throughout the paper where this hypothesis is 
fundamental. Although the construction can be done in the general 
case, it would fail to product a compactification in most of the cases. To avoid some technical complication, we will restrict ourselves
to moduli spaces over $\ZZ_p$ (so that we can focus on only one
prime). This restriction can be lifted but only at the
cost of even more technicalities in the definition, which did not 
seem worth the effort.
More precisely, let $g, g', p \in \NN$ with $p$ prime and $g \ge 2$ and 
denote by ${\mathscr H}^c_{g, g', p}$ the moduli space over $\ZZ_p$ of 
morphisms of degree $p$ between smooth curves of genus $g$ and $g'$
with constant ramification. That is, the degree $p$ geometrically
separable morphisms $f:X \to Y$ between proper smooth $S$-curves, 
where $S$ is a $\ZZ_p$-scheme, $g(X)=g$ and $g(Y)=g'$ and such 
that $\coker df$ is locally of the form $\O_S[x]/(x^n)$ for the finite 
flat topology. In particular, we have a canonical morphism
${\mathscr H}^c_{g, g', p} \to {\mathscr H}_{g, g', p}$ which is easily proved
to be bijective on the geometric points, but is not an isomorphism. 
Hence we can consider ${\mathscr H}^c_{g, g', p}$ as a stratification 
of ${\mathscr H}_{g, g', p}$.

It is relatively easy to see that ${\mathscr H}^c_{g, g', p}$ splits 
in two components: one that lives over $\ZZ_p$ and one 
that lives over $\FF_p$. More precisely, denote
by ${\mathscr H}^{c, \infty}_{g, g', p}$ the substack of 
${\mathscr H}^c_{g, g', p}$ parameterising morphisms $f:X \to Y$
over $S$ such that there exists a point $\p \in X$ at which
$f$ is of degree $p$ and $\coker df$ is, locally at $\p$ for the
finite flat topology on $X$, of the form $\O_S[x]/(x^n)$
with $n > p-1$. In particular, if $S$ is a field, a morphism $f$
is in ${\mathscr H}^c_{g, g', p}$ if and only if $S$ is of characteristic
$p$ and there exists a point at which $f$ is of degree $p$.
As we require the ramification locus to be constant, such
a morphism cannot be lifted to characteristic $0$. In particular,
this defines an open and closed substack of ${\mathscr H}^c_{g, g', p}$.
We will denote by ${\mathscr H}^{c, <\infty}_{g, g', p}$ the complement
of ${\mathscr H}^{c, \infty}_{g, g', p}$ in ${\mathscr H}^c_{g, g', p}$.

In this article, we construct a modular compactification of 
${\mathscr H}^c_{g, g', p}$ with only mild (and completely explicit)
singularities. This is achieved using a quite technical definition
of the admissible covering which generalises the existing ones.
The main ingredient is the concept of \emph{separating data}
which is a generalisation of the different in the global case
and over a general basis. Another important feature of the
admissible covering is the concept of the \emph{Hurwitz graph}
which was introduced originally by Henrio 
(cf. \cite{Henrio_p_adique}) but is vastly modified here to fit our 
purpose. The
Hurwitz graph is a combinatorial data which is attached on each
geometric fiber of the morphism and which we need to ensure
the existence of a smooth deformation. What makes the
definition quite tricky is that we need quite a large amount of
compatibility conditions between the separating data and the
Hurwitz graph. Roughly speaking, all those conditions
are necessary to ensure that the tangent space of the deformation 
functor is as little as possible, i.e. to ensure that
the singularities of the moduli space are as mild as possible.

\smallskip

Let's now describe the structure of the article. In the second
section, we introduce and study some morphisms between differential
modules. The original problem is that, for a morphism $f:X \to Y$
over a scheme $S$, the homomorphism $df:\Omega_{Y/S} \to f_* \Omega_{X/S}$
may be zero at some points. The $p^r$-earnest morphisms
defined  in this section are a mean to bypass this problem and to 
keep track of the differential information.

\smallskip

In the third section, we turn to the definition and the basic
properties of the admissible covers. The definition
is highly technical and requires many steps and many intermediary
objects. This definition was thought out to enjoy good degeneracy 
and deformation properties (in fact, it was really constructed
by looking precisely at the degeneracy of morphisms and keeping
the properties that can be deformed in a nice way). 
At the end of the section, we compare our notion with the already
existing notions of admissible covers.

\smallskip

The fourth section is entirely devoted to the reduction
of separable morphisms between smooth curves. More precisely,
we prove that admissible coverings are stable by degeneracy.
The first subsection of this section
consists of the study of morphisms between formal annuli;
the second one uses the results of the first one to prove
the main result of this section.

\smallskip

The fifth section is concerned with the deformation theory of 
admissible covers and the deformation functor $D_{adm}$. For this, 
we construct another functor $D_{abs}$ and a morphism
$abs:D_{adm} \to D_{abs}$. The results are then the following: 
the functor $D_{abs}$ has a versal deformation, which
is explicit and allows us to give a thorough description
of the singularities. Furthermore the morphism $abs$ admits a 
relative versal deformation, which will simultaneously be
proved to be smooth. The relative versality of the morphism
$abs$ is obtained by a good deformation theory and theorems
à la Schlessinger, adapted to the relative situation.

\smallskip 

In the sixth part, we define the moduli space of admissible covers,
prove its representability by an algebraic stack (using Artin's
theorem), the density of ${\mathscr H}^c_{g, g', p}$ and its properness.

\smallskip

Finally, in the appendix, we gather some results about various
deformation theories.

\medskip

\noindent{\bf Notation and convention:} The letter $p$ will denote a fixed prime 
number and $\nu_p$ will be the $p$-adic valuation on $\bar \ZZ_p$ 
normalized by $\nu_p(p)=1$.

All the schemes and rings in this article will be $\ZZ_p$-schemes
or $\ZZ_p$-algebras.

Let $\ell$ be an integer and $X$ be a scheme in which $p^\ell=0$.
We will denote by $F^\ell$ the absolute Frobenius defined by 
$x \mapsto x^{p^\ell}$. This morphism is a homeomorphism. In particular,
the functor $F^{\ell}_*$ is exact in the category of sheaves of abelian
groups.

When we speak of a semi-stable curve, we don't always mean it is 
proper.

Let $A$ be a complete local ring and $B$ the completion of
a singular double point over $A$. A set of coordinates of $B$ 
is a couple $(x, y) \in B^2$ such that $b:=xy \in A$ and the morphism
$A[[X, Y]]/(XY-b) \to B$ defined by $X \mapsto x$, $Y \mapsto y$ is an
isomorphism.

All the considered graphs will be directed. For any edge $e$ 
we will denote by $o(e)$ the origin and $t(e)$ the target.

Let $i , j \in \ZZ$, $i \not = 0$ and write $r:=\frac{j}{i}$. 
We will say that an element $\alpha \in \bar \ZZ_p$ is a $r$-th root
of $p$ if $\alpha^i-p^j=0$. Two $r$-th roots of $p$
differ from unit (in particular, for divisibility conditions, 
we can pick anyone). We will denote by $p^r$ one of these elements
and denote by $\ZZ_p[p^r]$ the $\ZZ_p$-sub-algebra of $\bar \ZZ_ü$ 
generated by $p^r$. By abuse of notation, we will denote by 
$p^\infty=0$ and $\ZZ_p[p^\infty]=\FF_p$.

\section{Sheaf of $p^r$-earnest morphism}

For this whole section, we will fix an $r \in \QQ_+ \cup \{\infty\}$ 
verifying $r\le 1$ or $r=\infty$. All considered schemes will be 
$\ZZ_p[p^r]$-schemes.
Let us fix a scheme $S$, a finite morphism $f:X \to Y$ 
of degree $p$ between smooth and \emph{not} necessarily proper 
curves. Our goal is to study a certain class of homomorphisms
which belongs to $\Hom_{\O_X}(f^*\Omega_{Y/S}, \Omega_{X/S})$
(philosophically, if $r \le \infty$ those are the homomorphisms of the 
form $\frac{df}{p^r}$ up to a lifting in characteristic $0$).

The first step of the definition comes from an explicit description
of the formal power series case. This definition could probably
be done in a more abstract (and perhaps canonical) way. The general 
definition follows from algebraisation and flat base change.

\begin{definition}[Formal power series]
~\\
Let $r \in (\QQ \cap [0, 1]) \cup \{\infty\}$ and $A$ be a $\ZZ_p[p^r]$-algebra,
$B=A[[x]]$, $C=A[[u]]$, $f:B \to C$ an
$A$-morphism and $\widehat\Omega_{B/A}$ (resp. $\widehat \Omega_{C/A}$) the 
completion (for the $(x)$-adic (resp. $(u)$-adic) topology) of 
the module of differentials. Write $f(x)=\sum_i a_i u^i$.
Let $\delta:\widehat \Omega_{B/A} \to \widehat \Omega_{C/A}$ be a 
$B$-homomorphism and write $\delta(dx)=\left(\sum_j b_{j+1} u^j\right) du$.
If $r< \infty$, we will say that $\delta$ is \emph{$p^r$-earnest} if for 
all $i \in \NN$ the following property is true
$$\left\{
\begin{array}{l}
\textrm{if } i \not = 0 \mod p \textrm{ then } a_i=\frac{p^r b_i}{i} \\
\textrm{if } i=0 \mod p \textrm{ then } b_i=\left(\frac{i}{p^r}\right)a_i.
\end{array}\right.$$

If $r=\infty$ (in particular $A$ is a $\FF_p$-algebra) and $g \in A$, 
we will say that $\delta$ is \emph{$p^\infty$-earnest for $g$} if
for all $i \not = 0 \mod p$ we have $a_i=\frac{g b_i}{i}$ and $b_i=0$ if $i=0 \mod p$.

We will denote by $\Xi^r_A$ the set of $p^r$-earnest homomorphisms.
\end{definition}

The following lemma explains the behaviour of $\Xi^r_A$ with respect
to base change and its link to the elements of the form 
$\frac{df}{p^r}$. In particular, it proves that this notion is
well defined (i.e. does not depend on the choices).

\begin{lemma}
\label{PropertiesCompletion}
With the notation of the definition, the following properties
are true.
\begin{enumerate}
\item Let $A \to A'$ be a faithfully flat morphism, then
we have an exact sequence
\begin{equation}
\label{FormalDescent}
0 \to \Xi^r_A \to \Xi^r_{A'} \rightrightarrows \Xi^r_{A' \otimes_A A'}.
\end{equation}
\item Suppose that $A$ is local noetherian and complete with residue
field of characteristic $p$ and that $r < \infty$. Let $W(A)$ be a 
complete local noetherian $\ZZ_p[p^r]$-algebra endowed with a 
surjective morphism $W(A) \to A$ such that $p^r$ is regular in $W(A)$
(such an algebra exists thanks to \cite{BourbakiAC2}, Th\'eor\`eme 
IX.5.3). Then an element 
$\delta \in \Hom_B(\widehat \Omega_{B/A}, \widehat \Omega_{C/A})$ 
is $p^r$-earnest if and only if there exists a lifting 
$\tilde f:W(A)[[x]] \to W(A)[[u]]$ of $f$ such that $p^r | d\tilde f$ and 
$\delta$ is the reduction of $\frac{d\tilde f}{p^r}$.
\item Suppose that $A$ is local noetherian and complete, of characteristic $p$ and that $r = \infty$. Let $g \in A$ and
let $W(A)$ be a complete local noetherian ring of characteristic $p$
endowed with a surjective morphism $W(A) \to A$ and a lifting $\tilde g$
of $g$ which is regular in $W(A)$. Then an element 
$\delta \in \Hom_B(\widehat \Omega_{B/A}, \widehat \Omega_{C/A})$ 
is $p^\infty$-earnest if and only if there exists a lifting 
$\tilde f:W(A)[[x]] \to W(A)[[u]]$ of $f$ such that $\tilde g | d\tilde f$ and 
$\delta$ is the reduction of $\frac{d\tilde f}{\tilde g}$.
\item If $A$ is local noetherian, then the property to
be $p^r$-exact is independent of the uniformising parameters $x$ and
$u$.
\end{enumerate}
\end{lemma}

\proof 1. Write $A''=A' \otimes_A A'$. By the faithful flatness of 
$A \to A'$ we have an exact sequence
\begin{multline*}
0 \to \Hom_B(\widehat \Omega_{B/A}, \widehat \Omega_{C/A}) \to \Hom_{B\otimes_A A'}(\widehat \Omega_{B\otimes_A A'/A'}, \widehat \Omega_{C\otimes_A A'/A'}) \\ 
\rightrightarrows \Hom_{B\otimes_A A''}(\widehat \Omega_{B\otimes_A A''/A''}, \widehat \Omega_{C\otimes_A A''/A''}).
\end{multline*}
(even though the rings and modules do not commute with base 
change).
In particular, we get the exactness on the left in 
\eqref{FormalDescent}. Then we only have to prove that
the divisibility conditions descend but this is a classical result
of faithful flatness.

2. The ``if'' part is trivial (this is the definition of 
$p^r$-earnestness). Let us prove the other direction.
For all $i \not = 0 \mod p$, let $\tilde b_i$ be a lifting of $b_i$ in
$W(A)$ and for all $i = 0 \mod p$ let $\tilde a_i$ be a lifting of $a_i$
in $W(A)$. Then define
$$\tilde f=\sum_{i \not = 0 \mod p} p^r \frac{b_i}{i} u^i + \sum_{i = 0 \mod p} \tilde a_i u^i,$$
this is a lifting of $f$.
Then we see that, by definition,
$$d\tilde f=p^r \left(\sum_{i\not =0 \mod p} b_i u^{i-1} + \sum_{i = 0 \mod p} \frac{i}{p^r} \tilde a_i\right) du$$
hence the desired property is satisfied.

3. The proof is similar to the second case.
 
4. Thanks to the first property, we can perform a base change and 
suppose that $A$ is complete. Then the result is obtained by the 
second and third assertions. \qed

Thanks to this lemma, we can extend the notion of $p^r$-earnestness
to points which are not necessarily rational in the following way.

\begin{definition}
Let $A$ be a complete local ring, $B \to C$ a morphism of local
complete $A$-algebras which are formally smooth and of relative 
dimension $1$. Let $\widehat \Omega_{B/A}$ and $\widehat \Omega_{C/A}$
be the completion of the differential modules, $r \in \QQ_+ \cup\{\infty\}$
and $\delta \in \Hom_B(\widehat \Omega_{B/A}, \widehat \Omega_{C/A})$.
Then $\delta$ is said to be $p^r$-earnest if there exists
a finite faithfully flat $A[p^r]$-algebra $A'$ so that the closed 
points of $C \otimes_A A'$ are rational and the image of $\delta$ in
$\Hom_{B \otimes_A A'}(\widehat \Omega_{B\otimes_A A'/A'}, \widehat \Omega_{C \otimes_A A'/A'})$ is $p^r$-earnest.
\end{definition}

The lemma \ref{PropertiesCompletion} ensures that this
notion is well defined.

We are now able to give a definition in the general situation.

\begin{definition}
Let $S$ be a scheme and $f:X \to Y$ a finite morphism of
degree $p$ between smooth curves. Let $\p \in X$ and 
$\delta \in \Hom_{\O_Y}(\Omega_{Y/S}, f_* \Omega_{X/S})$.
We will say that $\delta$ is $p^r$-earnest if it is $p^r$-earnest at
the completion at each point.
We will denote by $\Xi^r_{X/S}$ the sub-sheaf of 
$\hom_{\O_Y}(\Omega_{Y/S}, f_* \Omega_{X/S})$ composed 
by the homomorphisms which are $p^r$-exact at each point. 
\end{definition}

\begin{remark}
If $r=0$, we see that the set of $p^r$-earnest morphism is reduced to
$df$ and is therefore of no significant interest.

If $r > 0$ we see that $p^r | df$, hence the morphism $f$ should
be purely inseparable where $p=0$. Contrary to the case $r=0$, 
the $p^r$-earnest morphism need no longer be unique.
\end{remark}

Let us state some properties of $p^r$-earnest morphisms which are
formal consequences of the lemma \ref{PropertiesCompletion}.

\begin{proposition}
Let $S$ be a scheme and $f:X \to Y$ be a morphism of smooth $S$-curves.
Then the notion of $p^r$-earnestness is local for the finite flat
topology on $Y$ and the flat topology on $S$.
\end{proposition}

Next thing we would like to give a description of the
$p^r$-earnest homomorphisms in the global case. For that, it is 
convenient to work with the Frobenius (when available).

Let $p$ be a prime and $S$ be a scheme where $p^\ell=0$.
Then we can define a Frobenius morphism $F^\ell:S \to S$
by $\alpha \mapsto \alpha^{p^\ell}$. As in the case of a field,
we can define a relative Frobenius $F^\ell_S:X\to X^\ell$ on each 
$S$-scheme. We have the following lemma which gives us the 
structure of this morphism.

\begin{lemma}
\label{DecompositionRelativFrobenius}
Let $X \to S$ be a smooth curve. Then the morphism
$F^\ell_S : X \to X^\ell$ is locally (over $X^{\ell}$) of
the form $B \to B[u]/(u^{p^\ell}-v)$ where $dv$ is a basis of 
$\Omega_{B/\O_S}$.
\end{lemma}

\proof Locally on $X^\ell$, we can find $v \in \O_{X^\ell}$ which is
a basis of $\Omega_{X^\ell/S}$ and which is a $p^\ell$-th power 
in $\O_X$. The lemma follows easily from this remark. \qed

In particular, consider the case of a local artinian base
of residue characteristic $p$ and a finite morphism $f:X \to Y$ of 
degree $p$ between smooth curves which is purely inseparable
at the special fiber.
Then writing things explicitly, we see that the induced
morphism $X \to Y^{\ell}$ can be identified with $F^{\ell+1}$.
That is, for all such $f$ we have a commutative diagram
$$\xymatrix{
X \ar[rd]_{F^{\ell+1}} \ar[r]^f & Y \ar[d]^{F^\ell} \\
 & Y^\ell
}$$

We can now give the following ``global'' criterion for a 
homomorphism to be $p^r$-earnest.

\begin{lemma}
Let $S$ be a local scheme where $p^\ell=0$ and $f:X \to Y$ a finite
morphism of degree $p$ between $S$-smooth curves which is purely
inseparable at the special fiber. Suppose that 
$Y^{\ell}=\Spec B$ and that we can write $X=\Spec B[u]/(u^{p^{\ell+1}}-v)$ and 
$Y=\Spec B[x]/(x^{p^\ell}-v)$.
In particular, $dx$ is a basis of $\Omega_{Y/S}$.
Let $\delta \in \Hom_{\O_Y}(\Omega_{Y/S}, f_* \Omega_{X/S})$ and write 
$f(x)=\sum_{i=0}^{p^{\ell+1}-1} a_i u^i$ and 
$\delta(dx)=(\sum_{i=0}^{p^{\ell+1}-1} b_i u^i )du$.
If $r < \infty$, then $\delta$ is $p^r$-earnest if and only if
for all $i$ the following is true
$$\left\{
\begin{array}{l}
\textrm{if } i \not = 0 \mod p \textrm{ then } a_i=\frac{p^r b_i}{i} \\
\textrm{if } i=0 \mod p \textrm{ then } b_i=\frac{i}{p^r}a_i.
\end{array}\right.$$
If $r = \infty$ then $\delta$ is $p^\infty$-earnest if 
for all $i \not = 0 \mod p$ we have $a_i=\frac{g b_i}{i}$.
\end{lemma}

A decomposition of $X$ and $Y$ as required in the lemma always exists
locally on $Y^{\ell}$, cf. lemma \ref{DecompositionRelativFrobenius} 
and the discussion above.

\proof The lemma is trivially true where $x$ is a uniformising 
parameter.
Therefore it is true over all the points of $Y^\ell$ where $v$ is a 
parameter. Let $\p$ be a point of $Y^\ell$. By definition, the notion
of $p^r$-earnestness can be checked after and étale base change on 
$Y^\ell$ (hence we can suppose that $Y^\ell$ is an open subset of the 
affine line $\Spec A[x]$ where $x$ is the parameter above) and stable 
under finite base change of $S$, hence we can suppose that $\p$ is 
rational. Thus there is a parameter of $Y^\ell$ at $\p$ which
is of the form $x-\alpha$ with $\alpha \in A$. The result can be obtained
in a straightforward way. \qed

\begin{corollary}
\label{Lifting-pr-earnest}
Let $A' \to A$ be a surjective morphism of local artinian rings with
residue characteristic $p$. Let $f:X \to Y$ be a finite
morphism between smooth $A$-curves which is purely inseparable
at the special fiber and let $\delta \in Hom_{\O_Y}(\Omega_{Y/S}, f'_*\Omega_{X/S})$ be a $p^r$-earnest homomorphism. Then there exists locally
on $Y$ a lifting $f':X' \to Y'$ of $f$ and a lifting $\delta'$ of 
$\delta$ such that $\delta'$ is $p^r$-earnest.
\end{corollary}

\proof Use the explicit description given by the lemma above and then
lift explicitly as in the proof of the lemma 
\ref{PropertiesCompletion}. \qed

In the following, we will need a lemma concerning exact forms which
will allow us to classify the possible liftings of $p^r$-earnest
homomorphisms. Since the techniques involved in the proof are close
to the one used in this chapter, we state it here.

\begin{lemma}
\label{AlgebraiseExactForm}
Let $k$ be a field of characteristic $p$, $X$ a smooth curve
over $k$ and $\omega\in \Omega_{X/k}$. Then $\omega$ is exact if and
only if for all $\p \in X$, $\omega$ is exact in $\widehat \Omega_{X/S, \p}$.
\end{lemma}

\proof Denote by $\Omega_{X/S}^c$ the sheaf of locally exact 
differential forms and $F:X \to X$ the absolute Frobenius.
We are going to prove that $F_* \Omega_{X/S}^c$ is coherent
(then using the general properties of coherent modules, it is
easy to get the result).

Denote by $d:\O_X \to \Omega_{X/S}$ the universal derivation.
By definition, $\Omega_{X/S}^c$ is the image sheaf of $d$.
As $F_*$ is an exact functor on the category of abelian sheaves
($F$ is an homeomorphism) we see that $F_* \Omega_{X/S}^c$ is the
image of the map $F_* d$. But $F_* d$ append to be $\O_X$-linear.
Hence $F_* \Omega_{X/k}^c$ is coherent being the image of a homomorphism
between coherent modules. \qed

\section{Basic properties of $\dif$-morphisms}

The present section defines the notion of $\dif$-morphisms. 
Roughly speaking, those morphisms are morphisms 
between semi-stable curves endowed with an additional differential
structure. This allows one to deform these morphisms generically in
a unique way, although the considered morphism may be inseparable
when restricted to some irreducible component. More precisely,
the $\dif$-morphisms will be morphisms with additional data satisfying
some conditions ``formal locally''.  The following sections will 
be devoted to degeneracy and deformations of $\dif$-morphisms.

Before we can give the definition of $\dif$-morphism, we 
need to introduce some other concepts.

All the considered graphs will be directed. For a directed graph
$\Gamma$, we will denote
by $\Ar(\Gamma)$ the set of directed edges (positive and negative),
by $\Ar^+(\Gamma)$ the set of positive edges and by $\Som(\Gamma)$ the 
set of vertices. For any directed edge
$e$ we will denote by $\bar e$ the edge with the opposite direction.

\begin{definition}[Hurwitz graph]
\label{HurwitzGraph}
Let $\Gamma$ be a directed graph and $k$ a field. A $k$-Hurwitz
data on $\Gamma$ is a couple $(m, r)$ with
\begin{enumerate}[i)]
\item $m: \Ar(\Gamma) \to \ZZ$ (the conductor);
\item $r:\Som(\Gamma) \to \QQ_+ \cup \{\infty\}$ (the earnestness degree)
\end{enumerate}
satisfying the following properties
\begin{enumerate}
\item $\forall e \in \Ar(\Gamma)$ we have $m(\bar e)=-m(e)$;
\item \label{Compatibilityrm} $\forall e \in \Ar^+(\Gamma)$ we have $r(o(e)) \le r(t(e))$ with equality if $m(e)=0$;
\item \label{SmallDegree}
we have either $\Im r \subset [0, 1]$ or $\Im r \subset \{0, \infty\}$;
\item for all $e \in \Ar^+(\Gamma)$ one has $m(e) \ge 0$;
\item \label{mPrimeTop} if $k$ is of characteristic $p$ then $m(e)$ 
is either $0$
or prime to $p$ for all $e \in \Ar(\Gamma)$ ;
\item if $k$ is of characteristic $0$ then $m=0$ and $r=0$.
\end{enumerate}
\end{definition}

This definition is inspired by Henrio (cf. \cite{Henrio_p_adique})
a modification of which was used in \cite{Papier1}. In fact,
our definition defers hugely from the original (we lost 5 functions
out of 6 and introduced a new one) but we give it the same name 
because the purpose of those data remains the same (namely, to 
ensure the existence of a lifting).

The definition above will be useful to study the deformation problem
but will not be enough to ensure a lifting. For a lifting to exist,
we need an additional condition on the graph which is given through 
the following definition.

\begin{definition}[Reduced Hurwitz graph]
~\\
Let $(\Gamma, m, r)$ be a Hurwitz graph. Denote by $\Gamma_{red}$ the
graph obtained from $\Gamma$ after contracting the edges $e$
which verifies $m(e)=0$ and identifying the vertices which are
linked by such edges. The functions $m$ and $r$ then descend to 
functions on $\Ar(\Gamma_{red})$ and $\Som(\Gamma_{red})$ (through
the property \ref{Compatibilityrm} of the definition 
\ref{HurwitzGraph}) which we will denote by $m_{red}$ and $r_{red}$, and
$\Gamma_{red}$ inherits an orientation from $\Gamma$. It is easily seen 
that $(\Gamma_{red}, m_{red}, r_{red})$ is a Hurwitz graph. It will be 
called the \emph{reduced Hurwitz graph} attached to $(\Gamma, m, r)$.
By definition the reduced Hurwitz graph verifies $\forall e$ 
$m_{red}(e) \not = 0$. 

We call the reduced Hurwitz graph \emph{good} if the
orientation on $\Gamma_{red}$ defines an order on $\Som(\Gamma_{red})$
(saying that $t(e) > o(e)$ for any $e \in \Ar^+(\Gamma_{red})$).
\end{definition}

\begin{remark}
\label{GoodHurwitzFunction}
We can quite easily prove that a reduced Hurwitz graph is good
if and only if there exists a function $\ell:\Som(\Gamma_{red}) \to \NN$
that verifies $\ell(t(e)) > \ell(o(e))$ for any $e \in \Ar^+(\Gamma_{red})$
and $\ell(s)=0$ for $s$ minimal.
For the lifting property, the existence of such a function
will be fundamental.

\medskip

If $k$ is of characteristic $0$, then the reduced graph is trivial.
\end{remark}

The following definition is mainly a computational help, 
all the rest may be independent of it, but it will simplify some
of the later proofs.

\begin{definition}[Distended morphism]
Let $S$ be a scheme and $X\to S$, $Y\to S$ be $S$-semistable curves. 
Denote by $\omega_{X/S}$ and $\omega_{Y/S}$ the relative
dualising sheaves. In particular we have canonical inclusions
$\Omega_{X/S} \to \omega_{X/S}$ and $\Omega_{Y/S} \to \omega_{Y/S}$. A 
\emph{distended morphism} is a couple $(f, \partial f)$ where
$f : X \to Y$ is an $S$-morphism and $\partial f:\omega_{Y/S} \to f_* \omega_{X/S}$
is a homomorphism of $\O_Y$-module such that the following diagram
is commutative
$$\xymatrix{
\Omega_{Y/S} \ar[r]^{df} \ar@{^(->}[d] & f_* \Omega_{X/S} \ar@{^(->}[d] \\
\omega_{Y/S} \ar[r]^{\partial f} & f_*\omega_{X/S}.
}$$
Very often, the morphism $\partial f$ will also be denoted by $df$
(although $df$ may not determine $\partial f$ in a unique way).
\end{definition}

A morphism between smooth curves is automatically endowed with
a unique distended structure but this fail to be true when one
considers semistable curves.
For example, consider the morphism $f:k[x, y]/(xy) \to k[u, v]/(uv)$
defined by $x=u^n$ and $y=v^m$ with $n \not = m$ in $k$. The
morphism $df$ cannot be extended to a morphism between 
dualising sheaves.

Moreover, an existing distended structure may not be unique.

\medskip

In the definition of $\dif$-morphisms, we will need some
invertible sheaves on curves which are trivial on the smooth locus.
The best way to formalise this for our purpose is given by the 
following definition.

\begin{definition}[ftsl-invertible sheaf]
~\\
Let $S$ be a scheme, $X \to S$ a semistable curve and $\L$ an
invertible sheaf on $X$. Denote by $U$ the smooth locus of $X \to S$.
We will say that $\L$ is formally trivial on the smooth locus (ftsl 
for short) if for all point $s$ of $S$,
the invertible sheaf induced by $\L$ on the formal completion
of $U \times_S \Spec \O_{S, s}$ along $s$ is trivial. That is, if
we denote by $\m_s$ the maximal ideal of $\O_{S, s}$, there
exists an isomorphism
$$\phi_U: \limproj_n\left( \L|_{U \times_S \O_{S, s}/\m_s}\right) \to  \limproj_n\left( \O_U \otimes_S \O_{S, s}/\m_s\right)$$
Let us choose, for any singular point $\p \in X_s$, a trivialisation
$\phi_\p:\L_\p \otimes_{\O_{X, \p}} \widehat\O_{X, \p} \to \widehat\O_{X, \p}$.
The data $(\phi_U, \phi_\p)$ will be called a formal trivialisation
of $\L$ at $s$.
\end{definition}

Suppose that $S$ is local and that $U$ is affine. It's easy
to see that $\L$ is ftsl if and only if the restriction
of $\L_U$ to the fiber over $s$ is trivial. Indeed in this
case the deformations of $\L_s$ are classified by the sheaf 
$\H^1(U_s, \O_{U_s})$ which is zero because $U_s$ is affine.

\medskip

\begin{definition}[Unfolded separating data]
\label{UnfoldedSeparating}
~\\
Let $S$ be a scheme and $f:X \to Y$ a distended morphism between
$S$-semistable curves. An \emph{unfolded separating data} relative to $f$
is a triplet $(\L, \delta, g)$ where $\L$ is a ftsl 
invertible sheaf on $X$, $g:\L \to \O_X$ is a homomorphism
of $\O_X$-module and $\delta:\omega_{Y/S} \to f_*(\omega_{X/S} \otimes_{\O_X} \L)$
is a homomorphism of $\O_Y$-module, satisfying the following conditions
\begin{enumerate}[a)]
\item the diagram
$$\xymatrix{
\omega_{Y/S} \ar[rr]^{df} \ar[rd]_\delta & & f_* \omega_{X/S} \\
& f_*(\omega_{X/S} \otimes_{\O_X} \L) \ar[ur]_{f_*(Id \otimes g)}
}$$
is commutative;
\item \label{deltaInjective} $f^*\delta$ induces an injective morphism
in every fiber (which is equivalent to say that the cokernel
of $\delta$ is finite and flat over $S$);
\item \label{deltaIsom} $f^*\delta$ induces an isomorphism in a 
neighbourhood of all singular point;
\item \label{SmallConductor} $\coker f^*\delta$ is locally (for the 
finite flat topology) isomorphic to $\O_{S}[x]/(x^n)$  at a point $\p$. 
The integer $n$ is called the horizontal ramification degree at $\p$.
If $f$ is of degree prime to $p$ at $\p$ then we see
that $n=\deg_\p f-1$;
\item \label{gNonZero} the homomorphism $g$ is zero over no point 
(although it can be zero in restriction to some irreducible 
component, hence not injective).
\end{enumerate}

The unfolded data will be called \emph{finite} (resp. \emph{infinite}) 
if for any point $\p$ such that $f$ is of degree $p$ at $\p$, the 
horizontal ramification degree at $\p$ is $< p$ (resp. $\ge p$ and 
$\not = -1 \mod p$).
\end{definition}

The set of unfolded separating data is naturally endowed with
an action of $\O_{S}^*$ defined by 
$\alpha.(\L, g, \delta):=(\L, \alpha g, \alpha^{-1} \delta)$. The separating 
data will later be defined as an equivalence class of unfolded
separating data under this action.

\begin{remark}
The separating data are here to remember the information about
the ramification which is normally encoded in $df$ but can be 
trivial in the inseparable case. More precisely, the ``horizontal''
ramification is encoded in $\coker \delta$ and the ``vertical'' 
ramification is encoded in $\coker g$.
\end{remark}

\medskip

Let $S$ be a scheme, $X \to S$ a semistable curve and suppose that
for any geometric point $\bar s \to S$ we are given a structure
of a Hurwitz graph on the dual graph $\Gamma_{\bar s}$ of $X_{\bar s}$ 
(in particular, we are given an orientation on $\Gamma_{\bar s}$).
Then the oriented edges of $\Gamma_{\bar s}$ are naturally
in bijection with the branches at singular points of $X_{\bar s}$.
We can thus talk of the origin and target branches. When referring
to a set of coordinates $(u_o, u_t)$ at a singular point, $u_o$
will be a parameter of the origin branch and $u_t$ of the target 
branch. The singular points will be identified with the positive
edge.

\begin{definition}[Unfolded $\dif$-morphism]
\label{UnfoldedMorphism}
~\\
Let $S$ be the spectrum of a local ring, $X\to S$ and $Y \to S$
semi-stable curves such that the singular points of the special
fiber of $X$ are rational, $f:X \to Y$ a distended morphism, $(m, r)$ 
a Hurwitz structure on the dual graph $\Gamma$ of the special fiber 
of $X$ and $(\L, g, \delta)$ an unfolded separating data relative to $f$.
Denote with ``$\hat ~$'' the
data obtained by completion along the special fiber and $U$ the
smooth locus of $X \to S$. We will say that
$(f, m, r, \L, g, \delta)$ is an unfolded $\dif$-morphism
at the closed point of $S$ if there exists a formal trivialisation 
$(\phi_U, \phi_\p)$ of $\L$ and for any singular point $\p$ of $X$ there 
exists a set of coordinates $(u_o, u_t)$ of $X$ at $\p$ such that
\begin{enumerate}[i)]
\item \label{gScalaire} for any connected open formal subscheme $V$ 
of $\widehat U$ the homomorphism $g \circ \phi_U^{-1}|_V$ is given by an 
element of $A$ and if $r(V) < \infty$ then this scalar is in 
$p^{r(V)} A^*$;
\item \label{gPointDouble} for all singular point $\p$ corresponding
to a positive edge $e$, the morphism $g \circ \phi_\p^{-1}$ is given by 
an element of the form $\dif u_0^{m(e)}$ with $\dif \in A$;
\item \label{fPointDouble} for all singular point $\p$, the morphism
$f$ is given in a formal neighborhood of $\p$ by 
$x_o \mapsto u^n_o \alpha$, $x_t \mapsto u^n_t \alpha^{-1}$ where $(x_o, x_t)$ is 
a set of coordinates of $Y$ at $f(\p)$;
\item \label{TransitionFunction} for any singular point $\p$ 
corresponding to a positive edge $e$, the automorphism 
$\phi_\p \circ \phi_U^{-1}$ of $\Frac \widehat \O_{X, \p}$ is defined on the 
origin branch (resp. target branch) by the multiplication
by $u_o^{m(e)}$ (resp. $u_t^{m(\bar e)}$);
\item \label{DeltaExact} for any irreducible open subset $V$ of 
$\widehat U$ the homomorphism 
$\delta \circ f_*(Id \otimes \phi_U):\omega_{f(V)} \to f_* \omega_V$ is 
$p^{r(V)}$-horizontal.
\item \label{DeltaPointDouble} Let $\p$ be a double point
and $r_o$ and $r_t$ the degree of earnestness of $\delta$ on
the origin and terminal branches. Then $(u_0^{m(\p)}\delta)|_{u_0 \not = 0}$ 
is $p^{r_o}$-earnest and $(u_t^{-m(\p)}\delta)|_{u_t \not = 0}$ is 
$p^{r_t}$-earnest.
\end{enumerate}

The unfolded $\dif$-morphism $(f, m, r, \L, g, \delta)$ will be called 
finite if
\begin{enumerate}
\item $\Im r \subset \QQ_+$;
\item $(\L,g, \delta)$ is a finite unfolded separating data;
\end{enumerate}

It will be called infinite if
\begin{enumerate}
\item $S$ is an $\FF_p$-scheme;
\item $\Im r \subset \{0, \infty\}$;
\item $(\L,g, \delta)$ is an infinite unfolded separating data;
\end{enumerate}

Let $S$ be a scheme and $f:X \to Y$ a distended morphism between
semistable $S$-curves. For any geometric point $\bar s \to S$ 
let us a Hurwitz data $(m_{\bar s}, r_{\bar s})$ on the dual graph 
of $X_{\bar s}$ and an unfolded separating data $(\L, g, \delta)$ relative
to $f$. We will say that $(f, m_{\bar s}, r_{\bar s}, \L, g, \delta)$ is
an unfolded $\dif$-morphism if it's a $\dif$-morphism at each point 
after eventually a finite flat base change.
It's easily seen that the definition is independent of this base 
change.

The morphism $f$ will be called the underlying morphism of 
$(f, m_{\bar s}, r_{\bar s}, \L, g, \delta)$. Moreover, 
$(f, m_{\bar s}, r_{\bar s}, \L, g, \delta)$ will be called a structure
of unfolded $\dif$-morphism relative to $f$.
\end{definition}

\begin{definition}[$\dif$-morphism]
~\\
A \emph{$\dif$-morphism} is an equivalence class
of unfolded $\dif$-morphism between proper semistable
curves under the action of $\GG_m$ on unfolded separating data.
\end{definition}

The next step is to prove the compatibility of this definition
with those already existing in the literature. Let's do
it first for separable morphisms between smooth curves satisfying
an extra condition.

\begin{proposition}
Let $S$ be a scheme, $f:X \to Y$ a morphism between proper smooth 
curves such that for any point $s \in S$, the morphism $f:X_s \to Y_s$
is separable. Suppose that $\coker df$ satisfies condition 
\ref{SmallConductor}) of definition \ref{UnfoldedSeparating}.
Then $f$ admits a unique structure of $\dif$-morphism.
Moreover, if $S$ is of equal characteristic $p$ then this
$\dif$-morphism is infinite and if $S$ has no point of characteristic
$p$ then it is finite.
\end{proposition}

\proof The morphism $f$ is trivially distended (and admits a
unique structure of a distended morphism). Moreover, it is easily
seen that $(f, 0, 0, \O_X, 1, df)$ is a $\dif$-morphism
because, $f$ being separable in each fiber, the homomorphism
$df$ is injective in each fiber.

Let us suppose that $(f, m, r, \L, g, \delta)$ is another
$\dif$-morphism. By definition, we have $g\circ \delta=df$
and as $\delta$ and $df$ are injective, so is $g$. Moreover, using
formal trivialisation, we see that $g$ is an isomorphism in each
fiber. In particular, we can assume $\L=\O_X$. Then we can see
$g$ as a global automorphism of $\O_X$, that is an element of 
$\H^0(X, \O_X^*)$. But since $X$ is proper and smooth over $S$, we get
that $g$ is an element of $\H^0(S, \O_S^*)$. \qed

With almost the same proof, we can see that the moderated admissible
covering (cf. for example \cite{TheseWewers}) admits a unique 
structure of a $\dif$-morphism. Namely we have the following 
proposition.

\begin{proposition}
Let $S$ be a scheme, $f:X \to Y$ a moderated admissible covering
between proper semistable curves. Then there exists a unique
structure of $\dif$-morphism on $f$.
\end{proposition}

\proof The proof is similar to the case of smooth curves when 
taking into account the description of the admissible covering
in neighborhoods of singular points. Here we don't need the
extra condition on $\coker df$ because of the tameness of the
morphism. \qed

Then we can turn to the definition of admissible coverings which
are $\dif$-mor\-phisms with a minimality condition.

\begin{definition}[Admissible covering of degree $p$]
\label{AdmissibleCovering}
~\\
Let $S$ be a scheme, $X \to S$ and $Y \to S$ semistable curves.
Let $\underline f$ be a $\dif$-morphism between $X$ and $Y$ such that the
underlying morphism is finite of degree $p$. Then $\underline f$ will 
be called admissible if the following conditions are true
\begin{enumerate}[i)]
\item for any geometric point $\bar s \to S$ the reduced Hurwitz
graph of $(\Gamma_{\bar s}, m_{\bar s}, r_{\bar s})$ is good;
\item \label{MinimalityModel}
for any geometric point $\bar s \to S$ and any irreducible
component $V \in X_{\bar s}$ of genus $0$, $V$ meets the rest of 
$X_{\bar s}$ in at least $3-\# \supp (\coker \delta|_V)$ points.
\end{enumerate}
\end{definition}

\section{Degeneracy of unfolded $\dif$-morphism}

The goal of this section is to prove the following ``stable 
reduction'' theorem.

\begin{theorem}
\label{StableReduction}
Let $R$ be a discrete valuation ring of fraction field $K$
and $\underline f:X \to Y$ an admissible covering of degree $p$  
between smooth proper $K$-curves with $g(X) \ge 2$ and separable
underlying morphism. Then there exists 
a unique model of $f$ proper over $R$ which admits a unique 
structure of a $\dif$-morphism extending $\underline f$. 
Moreover, if $\underline f$ is finite (resp. infinite) then
this morphism is also finite (resp. infinite).
\end{theorem}

To prove this theorem, we first need to conduct a thorough study
of morphism between formal annuli. Then we will proceed with the
proof.

\subsection{Morphisms between formal annuli}

Let $R$ be a complete discrete valuation ring, $B:=R[[x, y]]/(xy-b)$,
$C:=R[[u, v]]/(uv-c)$ and $f:\Spec C \to \Spec B$ a finite morphism
of degree $n$ \'etale at the generic fiber.
We restrict ourselves to the case where the residue field
of $R$ is of characteristic $p>1$ and $n$ is a multiple of $p$
(otherwise the results are well known, see for example 
\cite{TheseWewers}). We will denote by $\varpi$ a uniformising 
parameter of $R$, $\Omega_B$ (resp. $\Omega_C$) the completion
of the differential modules of $\Spec B \to \Spec R$ (resp. 
$\Spec C \to \Spec R$) and $\omega_B$ (resp. $\omega_C$) the completion
of the canonical morphism. In particular, we have a natural injection
$\Omega_C \subset \omega_C$ and $\frac{du}{u}=-\frac{dv}{v}$ is a basis of
$\omega_C$.

Finally, for the special fiber we have $f(x) \in (u^n)$ or
$f(x) \in (v^n)$. We can assume that $f(x) \in (u^n)$ and thus that 
$f(y) \in (v^n)$.

The goal of this subsection is to prove the following proposition
which does not require that $n=p$.

\begin{proposition}
\label{MorphismAnnuli}
Suppose that $p | n$. Up to the multiplication of $x$ by 
an element of $R^*$ there exists $\alpha \in C^*$ such that
$f(x)=u^n \alpha$ and $f(y)=v^n \alpha^{-1}$. In particular, $df$ can
be extended to a homomorphism $\omega_B \to \omega_C$ in a unique way.

Moreover, there exists a $\dif \in R$ and $m \in \NN$ such that the 
homomorphism $df:\omega_B \to \omega_C$ is the composition of a morphism
$\delta:\omega_B \to \omega_C$ which induces an isomorphism
$\omega_B \otimes_B C \to  \omega_C$, and the multiplication by $\dif u^m$ or 
$\dif v^m$.

Finally, we have the following alternative 

(A) \ \ \ \begin{tabular}{l}
 - if $R$ is of equal characteristic $p$ then $m$ is prime to $p$; \\ 
 - if $R$ is of unequal characteristic then $\nu_p(m) < \nu_p(n)$. 
Moreover \\ 
$m=0$ if and only if $(\dif)=(n)$ in $R$. 
\end{tabular}
\end{proposition}

This proposition will allow us to prove the properties concerning 
double points (that is properties \ref{deltaIsom}) of definition 
\ref{UnfoldedSeparating}) and properties \ref{gPointDouble}), 
\ref{fPointDouble}), \ref{TransitionFunction}) and 
\ref{DeltaPointDouble}) of definition 
\ref{UnfoldedMorphism}). It also gives the definition and 
properties of the function $m$ of a Hurwitz graph if $\nu_p(n) \le 1$.

The rest of this subsection is devoted to the proof of the 
above proposition and is broken down in several steps.

First of all, using the norm homomorphism $\Frac(C)^* \to \Frac(D)^*$,
we can prove that $b=c^n$ up to an invertible element of $R$.
Thus, after eventually replacing $x$ by a multiple, we can assume
that $b=c^n$. In particular $b=0 \mod c$.

\noindent{$\rhd$ We can write $x=uP_u+cP_v$ with $P_u \in R[[u]]$ and $P_v \in R[[v]]$.}

As $uv=c$ in $C$, it is enough to prove that $x=uP_u \mod c$. 
We are going to prove this by induction on 
$q$. For $q=1$ the result is already known. 

Suppose it is true for $q$ and that $c=0\mod \varpi^{q+1}$. 

Write 
$x=u (\sum_{i \ge 0} a_i u^i)+\varpi^q (\sum_{j \ge 0} b_j v^j)$, 
$y=v (\sum_{j \ge 0} c_j v^j)+\varpi^q (\sum_{i \ge 0} d_i u^i)$ (such a 
decomposition of $y$ exists by symmetry).

As $0=b=c=uv=\varpi^{2q}= \mod \varpi^{q+1}$ we get
$$0=xy=\varpi^q u \sum_{i \ge 0} d_i u^i \sum_{i \ge 0} a_i u^i+\varpi^{q} v \sum_{j \ge 0} b_j v^j \sum_{j \ge 0} c_j v^j \mod \varpi^{q+1}.$$
In particular,  
$\varpi^q v \sum_{j \ge 0} b_j v^j \sum_{j \ge 0} c_j v^i=0$ in $R/(\varpi^{q+1})[[v]]$
but in this ring, $v \sum_{j \ge 0} c_j v^j$ is regular (because
$y \not =0 \mod \varpi$). Hence we get $\varpi^q \sum_{j \ge 0} b_j v^j=0 \mod \varpi^{q+1}$.

\medskip

\noindent{$\rhd$ The element $dx$ is a multiple of $du$ in $\Omega_C$.}

Write $x=uP_u+cP_v$ as above. Then we have 
$dx=\frac{\partial (uP_u)}{\partial u}du+c\frac{\partial P_v}{\partial v} dv$. 
But on the other hand, we have $cdv=uvdv=-v^2du$ in $\Omega_C$.
Hence we get $dx=\left(\frac{\partial (uP_u)}{\partial u}-v^2\frac{\partial P_v}{\partial v}\right)du$.

\medskip

Before continuing, we need a lemma which is probably well known
but will be used in several occasions.

\begin{lemma}
Let $P=\sum a_i u^i \in R[u]$ be a unitary polynomial.
Consider it to be a function on $\Spec C$ (via 
$R[u] \subset R[[uv]]/(uv-c)$) and assume that it has no zeros on 
$\Spec C \otimes_R K$. Then there exists $\alpha \in C^*$ such that 
$P=u^d \alpha$.
\end{lemma}

\proof Since the condition $u^d | P$ is stable under base change, it is enough to prove the result where $P$ splits. Let us then write
$P=\prod(u-\alpha_i)$ with $\alpha_i \in R$. The annulus $C$ can be 
described (rigidly) as $\{u | \ |c| < |u| < 1\}$. Since $P$ has no roots
in $\Spec C \otimes_R K$, any root of $P$ either verifies
$|\alpha_i| \le |c|$ in which case 
$(u-\alpha_i)=(u-c\frac{\alpha_i}{c})=u(1-v\frac{\alpha_i}{c})$ the element 
$(1-v\frac{\alpha_i}{c})$ being invertible in $C$,
or $|\alpha_i| = 1$ in which case $u-\alpha_i$ is invertible in $C$.
This proves the lemma. \qed

The Weierstrass preparation theorem of Henrio 
(cf. \cite{Henrio_galoisien}, Lemme 1.6) allows one to generalize
in the following direction.

\begin{corollary}
\label{ZeroGenericFiber}
Let $f \in C$ be non zero at the special fiber.
Suppose that $f$ has no zeros on the generic fiber. Then
$f$ is of the form $u^d \beta$ or $v^d \beta$ with $\beta \in C^*$.
\end{corollary}

\proof Direct consequence of the above lemma and 
\cite{Henrio_galoisien}, Lemme 1.6. \qed

\medskip

\noindent{$\rhd$ There exists $\alpha \in C^*$ such that $y=v^n\alpha$.}

Let us write $dx=hdu$, which is possible due to the previous step.
Using the fact that $xy=b$ and $uv=c$, we get
$$bdy=hc\left(\frac{y}{v}\right)^2 dv.$$
The quantity $\left(\frac{y}{v}\right)$ has a meaning in $C$ due to
the first step of the proof.

By hypothesis, the morphism $f$ is \'etale on the generic fiber,
hence $dy$ is a basis of $\Omega_C \otimes_R K$. In particular,
$\left(\frac{y}{v}\right)$ has no zero on the generic fiber.
By corollary \ref{ZeroGenericFiber}, we get that
$y$ is of the form $v^n \alpha$ for an $\alpha \in C^*$
Here we use that $y$ is of the form $v^n \alpha$ in the special 
fiber.

\medskip

\noindent{$\rhd$ The homomorphism $df$ extends in a unique way to 
$\omega_B \to \omega_C$ and there exists $\dif \in R$ and $m \in \NN$, 
$\beta \in C^*$ such that $\frac{dx}{x}=\dif u^m\beta \frac{du}{u}$ or 
$\frac{dx}{x}=\dif v^m \beta \frac{du}{u}$.}

Since $C$ is integral and $\omega_C \otimes_R K=\Omega_C \otimes_R K$
if the morphism $df$ extends to $\omega_B \to \omega_C$, this happends
in a unique way. The fact that $\frac{dy}{y}$ is a regular element 
of $\omega_C$ is a formal consequence of the previous step which tells
us that $y=v^n \alpha$.

Let us now prove the last assertion. We can write
$\frac{dy}{y}=h\frac{dv}{v}$. But by hypothesis, $\frac{dy}{y}$ is 
a basis of $\Omega_C \otimes_R K$ because $f$ is \'etale at the generic
fiber, hence $h$ has no zero in $\Spec C \otimes_R K$. The corollary
\ref{ZeroGenericFiber} yields the result.

\medskip

\noindent{$\rhd$ The alternative $(A)$ is true.}

Assume that $\frac{dy}{y}=\dif v^m \beta \frac{dv}{v}$ (cf. previous step).

We first consider the case of an equal characteristic ring and 
write $y=v^n \alpha$. By hypothesis $p | n$, hence we have
$\frac{dy}{y}=(d\alpha) \alpha^{-1}=\dif v^m \frac{dv}{v}$.
Therefore we see that $m-1$ is the vanishing order of 
$\frac{d\alpha}{\dif}$.
As $R$ is of equal characteristic $p$, $d\alpha$ has not term of 
order $n$ with $n = -1 \mod p$. Hence $m-1 \not = -1 \mod p$, which is
what we claimed.

Now we continue with the more complicated unequal characteristic 
case.

For this, we work directly in $R[[v, \frac{\varpi}{v}]]$ and write
$y=v^n\alpha$ and 
$\alpha=\sum_i a_i \left(\frac{\varpi}{v}\right)^i+\sum_j b_j v^j$.
As $\frac{dy}{y}=n\frac{dv}{v}+\frac{d\alpha}{\alpha}=\dif v^m \beta \frac{dv}{v}$
we get
\begin{equation}
\label{mUnequalChar}
\dif v^m \beta\alpha=n+\alpha^{-1}(-\sum_i ia_i \left(\frac{\varpi}{v}\right)^i+\sum_j jb_j v^j)
\end{equation}
As $\beta \alpha$ is invertible, we see that if $m=0$ then $(\dif)=(n)$ 
in $R$.

Suppose that $\nu_p(m) \ge \nu_p(n)$. Looking at the degree $0$ term
in \eqref{mUnequalChar} we get that $\dif | n$ in $R$. But by
definition, we know that $(\dif)=(mb_m)$ in $R$ hence $m | \dif$ and by 
hypothesis $n | m$ in $R$. In particular, we find that
$(n)=(\dif)$. But the order of $\frac{1}{\dif} \left(n+\alpha^{-1}(-\sum_i ia_i \left(\frac{\varpi}{v}\right)^i+\sum_j jb_j v^j)\right) \mod \varpi$ 
(which is supposed
to be $m$) is then equal to $0$. Thus we have $m=0$.

We have also proved that $m=0$ if $(\dif)=(p^{\nu_p(n)})$. \qed

\subsection{Proof of Theorem \ref{StableReduction}}

First we are going to choose the different data and prove their
uniqueness. Then we will prove that they satisfy the condition
of the definition.

\noindent{$\rhd$ Choice and uniqueness of a model of $f$}

Let us first recall a definition (cf. \cite{LiuStableReduction}): 
a semistable model of $f$ is a finite morphism $\tilde f:\X \to \Y$ 
between semistable proper $R$-curves the generic fiber of which is 
isomorphic to $f$. Such a model is called a stable model if it is 
minimal (for the domination relation) between all the semistable 
models. 

As $g(X) \ge 2$, $f$ admits (up to a finite extension of $K$) a 
stable model $f_0:\X_0 \to \Y_0$, cf. \cite{LiuStableReduction}, 
Corollary 4.6. Up to another finite extension of $K$,
we can assume the ramification locus to be composed of rational
points $e_i \in X(K)$. Denote by $\X_1$ the minimal semistable model 
obtained from $\X_0$ and unfolding the specialisation of the $e_i$
such that the $e_i$ specialize in the smooth locus.
For this model to be semistable, we may need another extension
of $R$. Then it is easily seen that there exists a model of
$\Y_1$ dominating $\Y_0$ and a finite morphism $f_1:\X_1 \to \Y_1$
which extends $f_1$.

Moreover, as this model is minimal, we see that it verifies the 
condition \ref{MinimalityModel} of the definition 
\ref{AdmissibleCovering}. Indeed, if $V \subset \X_0$ is an irreducible
component of the special fiber which does not satisfy this condition,
then $f(V)$ is an irreducible component of genus $0$ which intersects
the rest of $\Y_0$ in at most $2$ points. Hence $V$ and $f(V)$
can be simultaneously contracted and thus define a model of $f$,
contradicting the minimality of $f_1$.

Finally, if there exists an admissible covering $\underline f$
which is a model of $f$, then the underlying morphism must
be the morphism we just constructed because the condition
\ref{deltaIsom}) of the definition \ref{UnfoldedSeparating} ensures 
that the $e_i$ specialize in the smooth locus and condition
\ref{SmallConductor}) of the definition \ref{UnfoldedSeparating}
says that the specialisation of the $e_i$ in the special fiber are 
distinct.

By definition, for any singular point $\p \in \X_1$ the restriction
of the morphism $f_1$ to the generic fiber of $\Spec \O_{\X_1, \p}$
is étale; in particular the proposition \ref{MorphismAnnuli}
proves that $f_1$ is distended and that there exists a unique 
morphism $\omega_{\Y_1/S} \to f_{1 *} \omega_{\X_1, S}$ extending $df$.

\medskip 

\noindent{$\rhd$ Choice and uniqueness of the separating data}

Denote by $D$ the divisor of ramification on $X$ and by $\bar D$
its closure in $\X_1$. As the $e_i$ specialize in the smooth locus,
$\bar D$ is a Cartier divisor. Denote by 
$\L=\omega_{\X_1/S} ^{-1} \otimes f_1^* \omega_{\Y_1/S}(\bar D)$. Using the 
homomorphism $df:f^*\omega_{\Y_1/S} \to \omega_{\X_1/S}$, we get a natural
morphism $\O_{\X_1} \to \L^{-1}$ and thus $g:\L \to \O_{\X_1}$ which
is an isomorphism at the generic fiber.

In particular, we see that the Weil divisor associated to 
$g:\L \to \O_{\X_1}$ is precisely the vertical ramification of the
morphism $f_1$ (and only the vertical part).

By definition, we see that $df$ factorises through $\omega_{\X_1/S} \otimes \L \overset{g}{\to} \omega_{\X_1/S}$.

Suppose now that there exists an admissible covering $\underline f$ 
which is a model of $f$. We have already seen that the underlying
morphism is then $f_1$. Let $(\L', g', \delta')$ be a separating
data relative to $\underline f$. Then the fact that
$df_1=g' \circ \delta$, and $\delta'$ is injective in the special 
fiber (cf. definition \ref{UnfoldedSeparating} property 
\ref{deltaInjective}) and that $g'$ is an isomorphism in the generic 
fiber (which comes from the property \ref{gScalaire}) of the
definition \ref{UnfoldedMorphism}) implies that the Weil divisor 
associated to $g':\L' \to \O_{\X_1}$ is the vertical ramification of 
the morphism $f_1$. Hence we see that $g$ and $g'$ are isomorphic.
As $\X_1$ is a semi-stable proper curve on $R$, there exists
an $\alpha \in R^*=\H^0(\X_1, \O_{\X_1})$ such that
$\alpha.(\L, g, \delta)=(\L', g', \delta')$. This proves the uniqueness.

\medskip 

\noindent{$\rhd$ Choice and uniqueness of the Hurwitz graph}

For any geometric point $\bar s \to S$ and for any singular point
$\p \in \X_{1, \bar s}$, the proposition \ref{MorphismAnnuli}
gives a function $m$. With the notation of this proposition,
if $\frac{dy}{y}=v^m \dif \delta(\frac{dy}{y})$ and $e$ is the edge 
at $\p$ corresponding to the branch $v$ then $m(e)=m$ and
$m(\bar e)=-m$. The properties \ref{gPointDouble}) and 
\ref{TransitionFunction}) of definition \ref{UnfoldedMorphism} ensure
the uniqueness of $m$.

Moreover, we can define an orientation on the graph so that
it is compatible with $m$. This orientation is unique
only at the vertices where $m \not = 0$.

Let us now define the function $r$. If $R$ is of equal 
characteristic, and $V$ is an irreducible component of the smooth 
locus of $X$, we define $r(V)$ to be $0$ if $g|_V$ is invertible and 
$r(V)=\infty$ otherwise. In this case it is easily seen that $r$ is 
then uniquely determined.

Let us assume that $R$ is of unequal characteristic. Let $\eta$
be the generic point of an irreducible component $V$ of the special
fiber of $\X_1$. Denote by $\varpi$ a uniformising parameter of $R$.
As the special fiber of $\X_1$ is reduced, $\varpi$ is a uniformising 
parameter of the discrete valuation ring $\O_{\X_1, \eta}$. Thus, there
exists an element $\dif_V \in R$ (the different of the extension) 
which is a generating parameter of 
the ideal $\L_\eta \subset \O_{\X_1, \eta}$. Denote by $\nu_{\varpi}$ a
valuation on $R$. We define $r(V):=\frac{\nu_{\varpi}(\dif_V)}{\nu_{\varpi}(p)} \in \QQ_+$. The property \ref{gScalaire}) of definition
\ref{UnfoldedMorphism} ensures the uniqueness of $r$
because $\delta$ has to be injective and $g\delta=df$.

\medskip

Now we only need to check the properties required for 
$(f_1, m, r, \L, g, \delta)$ to be an unfolded $\dif$-morphism.

\noindent{$\rhd$ The data $(m, r)$ define a Hurwitz graph}

We only have to check properties \ref{SmallDegree}, 
\ref{Compatibilityrm} and \ref{mPrimeTop} of definition
\ref{HurwitzGraph}. Concerning property \ref{SmallDegree}, 
this comes from the fact that we are looking at morphisms of degree
of at most $p$ and that $r(V)$ is defined to be the different.

Property \ref{Compatibilityrm} comes from the explicit 
description of the data in a neighbourhood of a singular point
given by proposition \ref{MorphismAnnuli}.

The last property comes from the alternative $(A)$ of
proposition \ref{MorphismAnnuli} using the fact that the degree
of the morphism is $\le p$.

\medskip

\noindent{$\rhd$ The data $(\L, g, \delta)$ is an unfolded separating data.}

Let us prove first that $\L$ is ftsl. For that, we can 
restrict ourselves to an open smooth irreducible formal subscheme $W$
of the formal completion of $\X_1$. Let $\eta$ be its generic point.
Via the morphism $g$, we can look at $\L_\eta$ as an ideal of
$\O_{W, \eta}$. As we have seen before, denoting $\varpi$ a uniformising
parameter of $R$, $\varpi$ is a uniformising parameter or $\O_{W, \eta}$.
Thus there exists $\ell \in \NN$ such that $\L_\eta=\varpi^\ell. \O_{W, \eta}$.
It is then easy to prove that $\L_{W}=\varpi^\ell \O_W$ (via $g$)
because both Weil divisors are equal at each point of codimension $1$
(they are both trivial at the generic fiber) and $U$ is normal.

The fact that $g \circ \delta=df$ comes from the definition of these 
objects. Moreover, again from the definition, we see that
$\delta$ is injective.

The homomorphism $f_1^* \delta$ induces an isomorphism in a neighborhood
of each singular point due to the explicit description of
proposition \ref{MorphismAnnuli}. Moreover, as we have taken 
a model that splits the specialization of the ramification locus,
if follows the fact that $\coker \delta$ is of the form prescribed by 
the property \ref{SmallConductor} of definition \ref{UnfoldedSeparating}. The bound appearing in the definition comes from the bound
in the hypothesis of our morphism.

Let us prove now that $g$ is not zero at the special fiber. For that,
assume the contrary. Then we get in particular that $df=0$ at
the special fiber. As the morphism is of degree $p$, it follows
that $f$ is an homeomorphism. Thus $g(X)=g(Y)$ but the Hurwitz formula
tells us that $2g(X)-2-p(2g(Y)-2) \ge 0$, which means that
$(1-p)(2g(X)-2) \ge 0$. This is absurd because $p \ge 2$ and $g(X) \ge 2$.

\medskip

\noindent{$\rhd$ The data $(f_1, m, r, \L, g, \delta)$ is an unfolded $
\dif$-morphism.}

First let us take for $\phi_U$ the trivialisation of $\L$ on the 
smooth locus obtained above such that $g \circ \phi_U^{-1}$ is given
by an element of $R$. Up to a finite extension of $R$
it is easy to find trivialisation of $\L_\p$ for each singular point 
$\p$ such that the conditions \ref{TransitionFunction}) and 
\ref{gPointDouble}) of definition \ref{UnfoldedMorphism} are 
satisfied. Simply look at the definition of $(\L, g, \delta)$ in a 
neighborhood of $\p$ using the local description given by the
proposition \ref{MorphismAnnuli}. This procedure
may require a change of the set of coordinate, thus leading
to the choice of $u_o$ and $u_t$.

Property \ref{fPointDouble}) of definition \ref{UnfoldedMorphism}
follows directly from proposition \ref{MorphismAnnuli}.

Finally, let $V \subset \widehat \X_1$ be an open connected formal
subscheme and $\dif$ be a generator of $\L|_V\subset \O_V$ which is in 
$R$. Suppose that $R$ is of unequal characteristic, then
up to the multiplication by an invertible element of $R$ we can 
assume that $\dif=p^{r(V)}$.
We see that, with respect to the trivialisation $\phi_U$, we
have $\delta=\frac{df}{\dif}$ and this implies $\delta$ being
$p^{r(V)}$-earnest. 

If $R$ is of equal characteristic, the corresponding property
is proved in the same way.

\medskip{$\rhd$ The data $(f_1, m, r, \L, g, \delta)$ define an admissible
covering.}

We already saw that $(f_1, \L, \delta)$ satisfies the property 
\ref{MinimalityModel}) of definition \ref{AdmissibleCovering}.
It is thus enough to look at the other condition.

Thanks to remark \ref{GoodHurwitzFunction}, it is enough to 
construct an ``increasing'' (relatively to the orientation) function
on the graph of the special fiber. 

Choose a valuation $\nu_\varpi$ on $R$. Then for any irreducible 
component of the special fiber of $\X_1$ denote by $\dif_V \in R$ 
the element corresponding to $(g \circ \phi_U^{-1})|_V$. 
Then define $\ell(V):=\nu_{\varpi}(\dif_V)$.
Proposition \ref{MorphismAnnuli} shows that the function
$\ell$ has the desired property. This proposition
describes precisely how the different changes going from
one branch to the other.

\section{Deformation Theory of $\dif$-morphisms}

In this section, we want to prove the prorepresentability
of the functor of deformations of admissible covers. 
More precisely, let $k$ be a field and $\underline f_0$ be an
admissible covering. If $\underline f_0$ is an infinite covering, 
denote by $\A$ the category of artinian local $k$-algebras with 
residue field $k$, if $\underline f_0$ is finite, denote by $\A$ the 
category of artinian local rings whose residue field is $k$. Then 
we can define 
a functor $D_{adm}:\A \to Set$ by 
$$D_{adm}(A)=\{\textrm{admissible covers } \underline f \textrm{ over } A \textrm{ with special fiber } \underline f_0\}/isom.$$

The aim of this section is to prove the following theorem.

\begin{theorem}
\label{AdmProrep}
Suppose that $p \not = 2$, then the functor $D_{adm}$ has a universal 
deformation.
\end{theorem}

The proof of this will be decomposed into several steps which will
give more information about this functor. In particular,
we will get a description of its singularities.

The first step is to prove that there is no infinitesimal
automorphism. This is done in the following lemma.

\begin{lemma}
\label{AutomorphismCovering}
Suppose $p \not = 2$.
Let $A' \to A$ be a small extension of artinian local rings
with residue field $k$ and kernel $\af$ and $\underline f':X' \to Y'$ 
a deformation of $\underline f_0$. Then $\underline f'$ has no non trivial
automorphism as lifting of $\underline f' \otimes_{A'} A$.
\end{lemma}

\proof Let $\sigma$ be an automorphism of $\underline f'$.
Then $\sigma$ can be decomposed in $(\sigma_X, \sigma_Y, \sigma_\L)$
where $\sigma_X$ is an automorphism of $X$, 
$\sigma_Y$ is an automorphism of $Y$ and $\sigma_\L$ is an
automorphism of $\L$, satisfying compatibility conditions (in 
particular, $\sigma$ should fix $\delta$).
As $X_0$ is proper and semi-stable, $\sigma_\L$ can be thought
of as an element of $1+\af$. Using usual deformation
theory of local complete intersection morphisms, one sees
that $\sigma_X$ and $\sigma_Y$ are given by vector fields
$\chi_X$ and $\chi_Y$ defined on the special fiber.
Let $V$ be an irreducible component of $X$. If $V$ is of genus $\ge 1$,
we see that $\chi_X$ must be zero. This implies that $\chi_Y=0$.
If $f$ is inseparable this comes from the fact that $f(V)$ is also 
of genus $\ge 1$. If $f$ is separable this comes from the 
compatibility condition as expressed in theorem 
\ref{DeformationMorphism}.

Thus we can assume $V$ to be of genus $0$. Assume first
that $f|_V$ is of degree $< p$. As $\sigma$ must preserve
$\delta$, we see that $\chi_X$ must have a zero at each point
of the support of $\coker \delta_0$. Thus by the condition
on the number of points of this support, we get that $\chi_X=0$.
As before, we deduce from this that $\chi_Y=0$ because $f$
is separable.

Suppose now that $V$ is of genus $0$ and that $f|_V$ is of degree
$p$. Then two cases can occur: either $f|_V$ is separable or
inseparable. If $f|_V$ is inseparable, it is easy to see (for example
by the theorem \ref{DeformationMorphism}) that $\chi_X=0$.
Moreover, if $V$ meets the rest of $X_0$ in at least $3$ points,
then it is a classical result that $\chi_Y=0$. Thus we can suppose that
$V$ intersects the rest of $X_0$ in less than $3$ points. In 
particular, as $\underline f_0$ is admissible, the support of 
$\coker \delta_0$ is not empty.
Let us choose an isomorphism $V \to \PP^1_k$ such that
the point $0$ is a point where $V$ meets the rest of
$X_0$ and the point at infinity is not a point of the support
of $\coker \delta_0$. Write $V\setminus\{\infty\}=\Spec k[x]$ such that
$x=0$ is a point where $V$ meets the rest of $X_0$ and write
$\chi_Y=(a_0+a_1 x+a_2 x^2)\frac{\partial}{\partial x}$. As $X_0$ has a
singular point at $x=0$, $\chi_Y$ must be zero at $0$. Hence $a_0=0$.
Then writing that $\sigma_Y$ fixes $\delta$, one get
$a_1=2a_2=0$ (by writing things explicitly). That is, $\chi_Y=0$.

If $f|_V$ is separable of degree $p$, then the result follows 
by analysing the different cases. For example, suppose that
$V$ meets the rest of $X_0$ in only one point. Then by the Hurwitz
formula and the fact that $\coker \delta$ has its support in at least 
$2$ points one gets the existence of at least one point in
the support of $\coker \delta$ at which $f$ is of degree $< p$.
In particular, $\chi_X$ must have a zero at this point.
Suppose that $f$ is of degree $p$ at another point (otherwise the 
result follows trivially).
Using the compatibility of $\sigma$ and $\delta$ one gets
that $\chi_X=-\chi_Y$. But the compatibility with $f$ says that
$\chi_X=\chi_Y$ thus (as $p \not = 2$) one has $\chi_Y=\chi_X=0$.\qed

A problem that arises when studying this functor is the
need to speak about $p^{r(e)}$-earnest morphisms and it is
not easy to describe these morphisms if $p^{r(e)}$ is not in the 
considered ring.
To bypass this problem, we can proceed as follows.
Let $k$ be a field and $\underline f_0$ be a finite $\dif$-morphism
(we don't need $\underline f_0$ to be an admissible covering).
Let $r_0 \in \QQ_+$ such that all the $r(e)$ ($ e \in \Som(\Gamma)$) are 
integer multiples of $r_0$.
Denote by $\A_{r_0}$ the category of $\ZZ_p[p^{r_0}]$-algebras which
are local artinian and with residue field $k$.

If $\underline f_0$ is an infinite $\dif$-morphism, we define
$r_0=\infty$ and $\A_{r_0}$ to be the category of artinian $k$-algebras
with residue field $k$.

Then consider the functor $D_{\dif}:\A_{r_0} \to Set$ defined by 
$$D_{\dif}(A)=\{\dif\textrm{-morphisms } \underline f \textrm{ over } A \textrm{ with special fiber } \underline f_0\}/isom.$$

\begin{theorem}
The functor $D_{\dif}$ admits a formal versal deformation. 
Moreover, if $p \not = 2$ and  $\underline f_0$ is an admissible 
covering, then $D_\dif$ admits a universal deformation.
\end{theorem}

Theorem \ref{AdmProrep} then follows by descent theory and
the fact that $D_{\dif}$ has a universal deformation in the
admissible covering case, which follows from the lemma above.

To study $D_\dif$, we are going to define a functor
$D_{\dif-abs}:\A_{r_0} \to Set$ endowed with a natural morphism
$abs:D_{\dif} \to D_{\dif-abs}$. We then prove that the functor
$D_{\dif-abs}$ has a versal deformation by exhibiting it,
which allows us to precise its geometry. The last subsection
is devoted to the morphism $abs$. Namely we prove that it admits
a relative formal versal deformation using Schlessinger's criterion
adapted to the relative situation and deformation theory. A
step of the proof is the formal smoothness of $D_{\dif-abs}$.

\subsection{Versality of the functor $D_{\dif-abs}$}

We define at this point the functor $D_{\dif-abs}$. Roughly speaking,
it encodes the information about the ftsl sheaf
$\L$, the homomorphism $g$ and the singularities of $X$. 
If the source $X$ of the
$\dif$-morphism $\underline f$ is smooth, then we know
that $\L$ is trivial and $g$ is an isomorphism.
In particular, we can assume that $\L=\O_X$ and $g=Id$.
Thus from now on we will assume $X$ to be not smooth. In particular,
the scheme $U:=X \setminus Sing(X)$ is affine.

Before going into the definition, we need a result
which will enable us to eliminate $\L$. This lemma will
be a corollary of a formal patching result of which we need only a
particular case.

\begin{theorem}
Let $A$ be a noetherian complete local ring, $X \to \Spec A$
a proper semi-stable curve. Denote by $\mathfrak M$ the category of
invertible sheaves on $X$, $U$ the smooth locus of $X \to \Spec A$,
$\mathfrak M^\circ$ the category composed by the triples $(\L_U, (\L_\p)_{\p \in Sin(X)}, (\phi_\p)_{\p \in Sing(X)})$ where $\L_U$ is an invertible
$\O_U$-module, $\L_\p$ is an invertible $\widehat \O_{X, \p}$-module
and $\phi_\p:\left(\L_U \otimes_{\O_U} \Frac \widehat \O_{X, \p}\right) \to \left(\L_\p \otimes_{\widehat \O_{X, \p}} \Frac \widehat \O_{X, \p}\right)$ is an isomorphism of $\Frac \widehat\O_{X, \p}$-module (the isomorphism of such triples being defined in an obvious way).
Then the canonical morphism
$$\L \mapsto \left(\L|_U, (\widehat \L_\p), (\L|_U \otimes_{\O_U} \Frac \widehat \O_{X, \p} \overset{Id}{\to} \widehat \L_\p \otimes_{\widehat \O_{X, \p}} \Frac \widehat \O_{X, \p} )\right)$$
is an equivalence of category.
\end{theorem}

\proof The proof is essentially similar to the one in  
\cite{HarbaterStevenson} Theorem 1 but we need to remove the 
hypothesis that $R$ is a ring of formal power series. 

More precisely, denote by $i_U:U \to X$ and $i_\p : \p \to X$,
then the result relies on the exact sequence
$$0 \to \L \to i_{U *} \L|_U \oplus \bigoplus_\p i_{\p *} \widehat \L_\p \to \bigoplus_\p i_{\p *} \left(\widehat \L_\p \otimes_{\widehat \O_{X, \p}} \Frac \widehat \O_{X, \p}\right) \to 0$$
For a proof of the exactness of this sequence over a complete local
ring, see for example \cite{Papier1}, Lemmes 4.5 and 4.6.
\qed

Let $\underline f_1$ and $\underline f_2$ be two deformations of 
$\underline f_0$ over a base $A$ with isomorphic sources $X$. Denote
by $\L_1$ and $\L_2$ the underlying ftsl invertible
sheaves and chose a trivialisation as in definition
\ref{UnfoldedMorphism}. At each singular point $\p$ there 
exists sets of coordinates $(u_1, v_1)$ and $(u_2, v_2)$ such that the 
transition functions are of the form prescribed by property 
\ref{TransitionFunction}) of definition \ref{UnfoldedMorphism}
with the same $m$ since it depends only on the Hurwitz graph.
In particular, if we suppose that $k$ is algebraically closed,
there exists an automorphism of $\widehat \O_{X, \p}$ which sends
$u_1^m$ to $u_2^m$ and $v_1^m$ to $v_2^m$ because $m$ is 
prime to the characteristic of $k$ or equal to $0$, in 
which case we have nothing to do.
Then the formal patching theorem above shows that
$\L_1$ and $\L_2$ are isomorphic. Let's denote it by $\L$.

Let's chose, for each singular point $\p \in X$, a set
of coordinate $(u_{o, \p}, u_{t, \p})$ and denote by 
$\varpi_\p=u_{o, \p} u_{t, \p} \in A$ the thickness.
As $\L$ is ftsl, we have $\L|_U \cong \O_U$ and
the theorem above yields an identification
$$\Hom_{\O_X}(\L, \O_X) \overset{\sim}{\to} \left\{
\begin{array}{l}
(g_U, g_p) | g_U \in \O_U, g_\p \in \widehat \O_{X, \p} \textrm{ satisfying } \forall b \in \{o, t\} \\
\forall \p\in Sing(X)  \ g_U=g_p u_{b, \p}^{m(b(\p))} \textrm{ in } \widehat \O_{X, \p}[u_{b, \p}^{-1} ]
\end{array}\right\}.$$

Denote by $\Gamma$ the dual graph of the special fiber of $X$ and
let us take an element $(g_U, g_\p)$ of the right hand side set.
Then $g_U$ can naturally be decomposed into $(g_v)_{v \in \Som(\Gamma)}$.

Then we see that an element $(g_v, g_\p)$ defines an element of
$\Hom_{\O_X}(\L, \O_X)$ by the above identification
if and only if for all $\p \in \Ar^+(\Gamma)$ we have 
$g_{o(\p)} \varpi^{m(\p)}=g_{t(\p)}$ and $g_\p = g_{o(\p)} u_{o(\p)}^{m(\p)}$.

Thus we can define a bijection
$$\Hom_{\O_X}(\L, \O_X) \overset{\sim}{\to} \left\{
(g_v)_{v \in \Som(\Gamma)} | g_v \in \O_{v\cap U} | \forall \p \in \Ar^+(\Gamma)\  g_{o(\p)} \varpi^{m(\p)}=g_{t(p)}\right\}$$
i.e. we can ``forget'' the information about $g_\p$ without
actually forgetting anything.

Let us now define a functor $D_{\L-hom}:\A_{r_0} \to Set$.
For any $A \in \A_{r_0}$, the set $D_{\L-hom}(A)$ is composed by the 
isomorphism classes of  triples $(\U, (\X_\p), (g_v))$
where $\U$ is a deformation of $U:=X\setminus X_{sing}$, $\X_\p$
is a deformation of $\Spec \widehat \O_{X, \p}$ and $(g_v)$
is an element of the set $\Hom_{\O_X}(\L, \O_X)$ as described above
which is a deformation of the morphism $g$ defined by $\underline f$.

Let's define a subfunctor $D_{\L, g}$ of $D_{\L-hom}$.
If $X$ is smooth, then $D_{\L, g}$ is trivial.
If $X$ is not smooth then for all singular
point $\p \in X$, fix a ``universal'' thickness of the versal
deformation of $\Spec \widehat \O_{X, \p}$ (i.e. an isomorphism
between the versal deformation and $W(k)[[\varpi_\p]]$) and set
$D_{\L, g}(A)=(\U_0 \times \Spec A, X_\p, (g_v))$
with $g_v \in A$ for all $A$ satisfying the conditions

\begin{equation}
\left \{\begin{array}{l}
\textrm{if } r(v) < \infty \textrm{ then } g_v=p^{r(v)} \\
\forall \p \in \Ar^+(\Gamma)\  g_{o(\p)} \varpi_\p^{m(\p)}=g_{t(p)}.
\end{array} \right .
\end{equation}

Then, by definition of an unfolded $\dif$-morphism
the canonical morphism $D_{\dif} \to D_{\L-hom}$ factorises through
$D_{\L, g} \to D_{\L-hom}$ thus giving rise to a morphism
$D_{\dif} \to D_{\L, g}$.

\begin{proposition}
\label{DLgDeform}
The functor $D_{\L, g}$ admits a versal deformation.
\end{proposition}

\proof This is clear if $X$ is smooth. If $X$ is not
smooth (in which case $U$ is affine) then all the deformations
of $U$ are trivial. Hence it's easily seen that a versal deformation
of $D_{\L, g}$ is given by the completion at ``$g_v$'' of
$$W(k)[p^{r_0}][[(\varpi_\p)]][(g_v)]/I,$$
where $I$ is given by the conditions above. \qed

Moreover, we see that the singularities are explicit and depend
only on the geometry of the \emph{reduced} Hurwitz graph.

Finally let us choose, for any singular point $\p$, a thickness 
$\varpi_\p$ of the universal deformation of $\Spec \widehat \O_{X, \p}$ and
denote by $n_\p$ the degree of $f_0$ at $\p$. Then
define a subfunctor $D_{abs}$ of $D_{\L, g}$ by the condition that
$\forall \p, \p'$ such that $f_0(\p)=f_0(\p')$ one has 
$\varpi_\p^{n_\p}=\varpi_{\p'}^{n_{\p'}}$

Then we immediately see that $D_\dif \to D_{\L, g}$ factorises through
$D_{abs}$ giving rise to a morphism $abs:D_{\dif} \to D_{abs}$, and
that $D_{abs}$ has a versal deformation.

\subsection{Versality of the morphism $abs$}

The aim of this subsection is to prove that the natural map $abs$
defined in the preceding subsection admits a relative formal versal 
deformation. The principal tool here is Schlessinger's deformation 
theory (cf. \cite{Deformation}) Although the original article
does not deal with the relative case but only the absolute, the
definitions and results can be generalised by the same methods.

It's therefore sufficient to prove the following theorem.

\begin{theorem}
\label{DeformationAbs}
Let $k$ be an algebraically closed field, $\underline f_0$ a 
$\dif$-morphism. Then there exists finite dimensional $k$-vector
spaces $T^0$ and $T^1$ such that for all $A \in \A_{r_0}$, any small
extension $A' \to A$ in $\A_{r_0}$ of kernel $\af$, any deformation
$\underline f$ of $\underline f_0$ over $A$ and any 
$\mathfrak g' \in D_{\L, g}(A')$ whose image in $D_{\L, g}(A)$ is 
$abs(\underline f)$, the following properties are true 
\begin{enumerate}[i)]
\item there exists a lifting $\underline f'$ of $\underline f$ over $A'$
such that $abs(\underline f')=\mathfrak g'$;
\item the set of isomorphism classes of liftings of $\underline f$ to
$A'$ with image $\mathfrak g'$ is endowed with a canonical action of
$\af \otimes_k T^1$ making it a principal homogeneous space;
\item the group of isomorphisms of a lifting of $\underline f$ to $A'$
is canonically isomorphic to $\af \otimes T^0$.
\end{enumerate}
In particular, by the first property, $abs$ is formally
smooth.
\end{theorem}

We are going to prove this theorem without the assumption
that the underlying curves are proper. Thus we will need to remove,
at first, the finiteness of the vector spaces. The
condition on the image of the deformation by $abs$ still
makes sense since it is equivalent to conditions on $\L$, $g$ and 
conditions on the deformation of the double points.

Let us first take care of the automorphism group of a lifting
(third part of the theorem above). For any invertible sheaf
$\F_0$ on $X_0$ and any section $s_0 \in \H^0(X_0, \F_0)$
denote by
$$T^0_{X_0, \F_0, s_0}:=\ker \left(\Ext^0_{\O_X}(\Pp^1_{X_0/k}(\F_0), \F_0) \overset{d^1_{s_0}}{\to} \H^0(X_0, \F_0))\right)$$
This sheaf classifies the automorphism of a deformation of 
$(X_0, \F_0, s_0)$ (see proposition \ref{AutLiftingSection}).
In particular, we have a canonical morphism, the forgetful functor
$(X_0, \F_0, s_0) \mapsto X_0$,
$$T^0_{X_0, \F_0, s_0} \to \Ext^0_{\O_{X_0}}(\Omega_{X_0/k}, \O_{X_0}).$$
Define the space $T^0_{X_0, \L_0, g_0, \delta_0}$ by the following 
Cartesian diagram :
$$\xymatrix{
T^0_{X_0, \L_0, g_0, \delta_0} \ar[r] \ar[d] & T^0_{X_0, \L_0^{-1}, g_0^\vee} \ar[d] \\
T^0_{X_0, f^* \omega_{Y_0(k)}^{-1} \otimes \omega_{X_0/k} \otimes \L_0, \delta_0} \ar[r] & \Ext^0_{\O_{X_0}}(\Omega_{X_0/k}, \O_{X_0}).
}$$
We see that it classifies the automorphism of a lifting of 
$(X_0, \L_0, g_0, \delta_0)$.
Moreover, we have a canonical morphism, the forgetful morphism
$(X_0, \L, g_0, \delta_0) \mapsto X_0$, 
$$T^0_{X_0, \L_0, g_0, \delta_0} \to \Ext^0_{\O_{X_0}}(\Omega_{X_0/k}, \O_{X_0}).$$

Then define $T^0_{\underline f_0}$ by the following Cartesian diagram
$$\xymatrix{
T^0_{\underline f_0} \ar[r] \ar[d] & T^0_{X_0, \L_0, g_0, \delta_0} \times \Ext^0_{\O_{Y_0}}( \Omega_{Y_0/k}, \O_{Y_0}) \ar[d]\\
\Ext^0(df_0, f_0) \ar[r] & \Ext^0_{\O_{X_0}}(\Omega_{X_0/k}, \O_{X_0}) \times_k \Ext^0_{\O_{Y_0}}( \Omega_{Y_0/k}, \O_{Y_0}).
}$$
Theorem \ref{DeformationMorphism} explains the bottom row.
It is then easily seen by construction that $T^0_{\underline f_0}$ 
satisfies the desired property. Moreover, through the explicit 
description above we see that if $X_0$ and $Y_0$ are proper, then
$T^0_{\underline f_0}$ is finite dimensional.

We will denote by ${\mathscr T}_{\underline f_0}$ the sheaf on $Y_0$ defined
for all $U$ by $T^0_{{\underline f_0}}|_{f^{-1}(U)}$.

We turn now to the two other parts of the theorem. To prove those,
we're going to prove that in each of the following cases, there 
exists a deformation which is unique up to isomorphism:
\begin{enumerate}
\item the affine smooth $U \subset X_0$ such that $f_0|_U$ is 
separable;
\item the affine smooth $U \subset X_0$ such that $f_0|_U$ is 
purely inseparable (as $f_0$ is of degree $p$, only those two cases
can occur in the smooth case);
\item the singular locus.
\end{enumerate}

Then as usual the global result follows defining 
$T^1_{\underline f_0}=\H^1(Y_0, {\mathscr T}_{\underline f_0})$.

The first point above 
is classical and can easily be deduced from the usual deformation 
theory of morphisms as stated in theorem \ref{DeformationMorphism}.
Hence we'll focus on the other two cases.

\subsubsection{Deformation of the inseparable smooth locus}

First, by corollary \ref{Lifting-pr-earnest} 
there exists a lifting (at least locally on $Y$).
Here we assume  $X_0$ and $Y_0$ to be smooth, affine,
irreducible and that $\Omega_{Y_0/k}$ admits a basis which is exact.
Let $\underline f'_1$ and $\underline f'_2$ be two liftings of
$\underline f$ with the same image under $abs$. We are going to prove 
that they are isomorphic. 
Considering the deformation theory of smooth affine pointed curves
(which is trivial) we see that we can set $Y':=Y'_1=Y'_2$, 
$X':=X'_1=X'_2$ and that $\coker \delta_1 =\coker \delta_2$.

Consider $g'_1$ and $g'_2$ as elements of $A'$;we can do so because
$X_0$ is smooth and irreducible. By hypothesis, 
$abs(\underline f'_1)=abs(\underline f'_2)$ hence $g':=g'_1=g'_2$. Moreover,
as $f_0$ is inseparable the element $g'$ is in the maximal ideal of 
$A'$. Hence $df'_1=g'\delta_1=g'\delta_2=df'_2$. Then proposition 
\ref{LiftingFrobUpToDiff} tells us that $f'_1=f'_2$ up to an 
isomorphism of $Y'$.

What is left to prove is that $\delta_1=\delta_2$. This
will be up to an isomorphism of $X'$. Denote by $\epsilon$ a 
generator of the ideal $\af$ and $r$ the exactness degree 
of $\delta_0$. First, it is convenient to recall how an infinitesimal
automorphism of $X'$ acts on $\delta_1$. For this,
suppose that $\Omega_{X'/A'}$ admits a basis of the form $du$
and denote by $\frac{\partial}{\partial u}$ the dual basis of $du$.
Then any infinitesimal automorphism $\chi$ is given by a vector
field of the form $\chi=\epsilon \alpha \frac{\partial}{\partial u}$.
Suppose moreover that $\Omega_{Y'/A'}$ admits a basis of the form
$dx$ and write $\delta_1(dx)=\beta du$.
Then it is easy to see that 
$$(\chi \delta_1)(dx)=\chi(\beta du)=\delta_1(dx)+\epsilon d(\alpha \beta).$$

As $\coker \delta_1=\coker \delta_2$ we see that
$\delta_1(dx)=\gamma \delta_2(dx)$ because they have the same zeros
with $\gamma=1+\epsilon \zeta$ (because $\delta_1$ and $\delta_2$ are 
liftings of $\delta$).
Therefore we get $\delta_1(dx)=\delta_2(dx)+\epsilon \zeta \delta_0(dx)$.
As $\delta_1$ and $\delta_2$ are earnest $\beta \delta_0(dx)$ is at least 
locally exact; to see this, prove it after completion and 
algebraise with the lemma \ref{AlgebraiseExactForm}. Thus we can
write $\beta \delta_0(dx)$ as $dh$. Moreover,
after eventually shrinking $X'$, we can choose $h$ such that
for all $\p \in X'$ one has $\ord_\p h \ge \ord_\p \delta_0(dx)$.
Denote by $\beta_0$ the reduction of $\beta$ modulo the maximal
ideal of $A'$ so that $\delta_0(dx)=\beta_0 du$ and write 
$\alpha:=\frac{h}{\beta_0}$. Then $\alpha$ is regular and 
the automorphism defined by $\epsilon \alpha \frac{\partial}{\partial u}$
sends $\delta_1(dx)$ to $\delta_2(dx)$. As $dx$ is a basis, it sends
$\delta_1$ to $\delta_2$.

\subsubsection{Deformation of the singular locus}

Through descent theory, we are reduced to the case of morphisms
$A[[x, y]]/(xy-a) \to A[[u, v]]/(uv-b)$ and by hypothesis
this morphism is of the form $x \mapsto u^d \alpha$, 
$y \mapsto v^d \alpha^{-1}$ with $\alpha \in  A[[u, v]]/(uv-b)$ invertible.
Denote by $\beta \in  A[[u, v]]/(uv-b)$ the invertible such that
$\beta \frac{du}{u}=\delta\left(\frac{dx}{x}\right)$. We'll suppose, to fix 
the ideas, that $g$ is of the form $\dif u^m$. Since we are looking 
at elements with fixed image under $abs$, the data $\dif'$ and the 
lifting $b'$ of $b$ (the thickness of the lifting) are fixed.
Then we are looking for the liftings of $\alpha$ and $\beta$
satisfying $d=g\delta$. Denote by $\gamma:=\alpha \beta$, and by 
$r_o$ and $r_t$ the degree of earnestness of $\delta$ on the
origin (say $(x)$) and the terminal (say $(y)$) branches.

Let us write $\alpha=\sum_{i \ge 0} \zeta_i u^i+\sum_{j > 0} \eta_j v^j$.
Then it is easy to see through an explicit computation that
$d=g\delta$ if and only if 
$$\sum_{i \ge 0} \zeta_i (p+i) u^i + \sum_{j > 0} \eta_j (p-j) v^j=\dif u^m \gamma.$$

Using the fact that $(u^m\delta)|_{x \not = 0}$ is $p^{r_o}$-earnest and 
$(v^{-m}\delta)|_{y \not = 0}$ is $p^{r_t}$-earnest, it is easy to see that
there exists a lifting of the morphism and of $\delta$ which
are compatible.

Although having to consider several cases, uniqueness of the lifting
can then be proved similar to the uniqueness in the smooth case.

\section{Properties of the moduli space of admissible covering}

In this section, we finally introduce the moduli space
of admissible covers, prove that it is representable by an algebraic
stack and that it is proper.

Let $g, g', p \in \NN$ with $g \ge 2$ and $p$ prime be distinct from $2$.
Denote by $\bar {\mathscr H}^c_{g, g', p}$ the groupoid over $\ZZ_p$ 
classifying the admissible covers of degree $p$ between 
curves of genus $g$ and $g'$. By the usual descent theory,
it is easy to see that $\bar {\mathscr H}^c_{g, g', p}$ is actually
a $\ZZ_p$-stack. 

Our main theorem is the following

\begin{theorem}
The stack $\bar {\mathscr H}^c_{g, g', p}$ is an algebraic Deligne-Mumford 
stack, contains ${\mathscr H}^c_{g, g', p}$ as an open dense substack, and
is proper over $\ZZ_p$.
\end{theorem}

\proof The fact that the diagonal morphism
$$\bar {\mathscr H}^c_{g, g', p}\overset{\Delta}{\to} \bar {\mathscr H}^c_{g, g', p} \times_{\ZZ_p} \bar {\mathscr H}^c_{g, g', p}$$
is representable, separated and quasi-compact follows from
the usual properties of $Isom$ schemes of curves.

Moreover we can prove by usual techniques that 
$\bar {\mathscr H}^c_{g, g', p}$ is locally of finite presentation. 
Then theorem \ref{AdmProrep}
together with Artin's algebraisation theorem (see for example
\cite{Champs}, Corollaire 10.11) shows that 
$\bar {\mathscr H}^c_{g, g', p}$ is an algebraic stack.
Moreover, by the fact that an admissible covering has no 
infinitesimal automorphism (cf. lemma \ref{AutomorphismCovering} and
the general deformation theory) one gets that the diagonal
morphism $\Delta$ above is non ramified and
$\bar {\mathscr H}^c_{g, g', p}$ is Deligne-Mumford.

\medskip
It is easy to see that ${\mathscr H}^c_{g, g', p}$ is open in
$\bar {\mathscr H}^c_{g, g', p}$. We now prove its density.
Let $k$ be a field and $\underline f_0$ an admissible
covering over $k$. We have to prove that $\underline f_0$ can
be deformed to a separable morphism between smooth curves.
In fact, by definition of an unfolded separating data
which tells us that the different $g$ is zero over no point,
it is enough to prove that $\underline f_0$ can be deformed
into a morphism between smooth curves. This morphism
will automatically be separable. The problem for $k$ of
characteristic $0$ is easy to solve, we assume $k$ to be 
of characteristic $p$.

By the formal smoothness of the functor $abs$ (cf. theorem
\ref{DeformationAbs}), it is enough to see that
there exists a lifting of $abs(\underline f_0)$ over a complete
discrete valuation ring $R$ such that the thickness at each point 
is non zero. If $\underline f_0$ is finite, then choose a 
ring $R$ of unequal characteristic with residue field $k$, if $R$ is 
infinite then choose $R=k[[t]]$.
As there is, up to an extension of $R$ no problem to deform the 
tame part of $abs(\underline f_0)$ (that is, the singular points 
where $f$ is of degree $< p$) we can focus on the wild part. As
stated after proposition \ref{DLgDeform}, the problem
depends only on the reduced Hurwitz graph, which we know is
good by definition of an admissible cover.
Suppose first that $\underline f_0$ is infinite (i.e. $R$ is of equal
characteristic).
In particular, due to remark \ref{GoodHurwitzFunction}, there
exists a function $\ell:\Som(\Gamma_{red}) \to \NN$ which verifies 
$\ell(t(e)) > \ell(o(e))$. Moreover, we can suppose that $\ell$ is
zero at the minimal points. Those are precisely the points
where $g$ is invertible. Let $\pi \in R$ be an element
with positive valuation.
Then define, for all $v \in  \Som(\Gamma_{red})$, $g_v:=\pi^{\ell(v)}$.
Denote by $m_0$ the gcd of all the $m(e)$ for $e$ an edge
of the reduced Hurwitz graph. Then, up to an extension of $R$,
we can suppose that there exists $\pi' \in R$ such that 
$\pi'^{m_0}=\pi$. Then for any edge $e$ of the reduced Hurwitz
graph, chose $\pi'^{\frac{m_0}{m(e)}}$ for the thickness.
It is then easily seen that it defines a lifting of 
$abs(\underline f_0)$.

In the finite case, repeat the procedure with $\ell=r$.

\medskip

What is left to prove is that $\bar {\mathscr H}^c_{g, g', p}$ is proper
over $\ZZ_p$. For this, we want to use the valuative criterion 
of properness. As in \cite{Deligne_Mumford} (remark after the
theorem 4.19, see also \cite{Champs}, Remarque 7.12.4)
it is enough to prove that for any discrete valuation ring
$R$ with quotient field $K$ and any morphism 
$\Spec K \to {\mathscr H}^c_{g, g', p}$, this morphism can be prolonged
up to a finite extension of $K$ in a unique way
to a morphism $\Spec R \to\bar {\mathscr H}^c_{g, g', p}$.
This version of the vlauative criterion for properness
is valid here because we have just proved that ${\mathscr H}^c_{g, g', p}$
is dense in $\bar {\mathscr H}^c_{g, g', p}$.

With this criterion in hand, the properness of 
$\bar {\mathscr H}^c_{g, g', p}$ follows from theorem 
\ref{StableReduction}. \qed

As in the case of ${\mathscr H}^c_{g, g', p}$, we can decompose
the algebraic stack $\bar {\mathscr H}^c_{g, g', p}$ in two open
and closed components : one containing the finite morphism
$\bar {\mathscr H}^{c, <\infty}_{g, g', p}$ and one containing
the infinite part $\bar {\mathscr H}^{c, \infty}_{g, g', p}$ (this last
part lives only in characteristic $p$).

\appendix
\section{Some results about classical deformations}

Let us first state a generalisation of a result of Ziv Ran 
(cf. \cite{ZivRan}) concerning the deformation theory of morphisms
between pointed semistable curves. 

\begin{theorem}
\label{DeformationMorphism}
Let $k$ be a field and $f_0:X_0 \to Y_0$ a morphism between semi\-stable
(not necessarily proper) curves. Let $A$ and $A'$ be a local 
artinian rings with residue field $k$, $f:X \to Y$ a deformation 
of $f_0$ and $A' \to A$ a small extension with kernel $\af$. Then
there exist $k$-vector fields $\Ext^i(df_0, f_0)$ fitting in a 
long exact sequence 
\begin{multline}
\label{ExactSequenceRan}
\ldots \to \Ext^{i-1}(\Omega_{Y_0/k}, f_* \O_{X_0})\to \Ext^i(df_0, f_0) \to \\
\Ext^i(\Omega_{Y_0/k}, \O_{Y_0}) \oplus \Ext^i(\Omega_{X_0/k}, \O_{X_0}) \to 
\Ext^i(\Omega_{Y_0/k}, f_* \O_{X_0}) \to \ldots
\end{multline}
such that
\begin{enumerate}
\item there exists a canonical element $\omega \in \af \otimes_k \Ext^2(df_0, f_0)$, the vanishing of which is necessary and sufficient for a
lifting of $f$ to $A'$ to exist;
\item if there exists a lifting of $f$ to $A'$, then the set of
isomorphism classes of liftings is naturally a principal homogeneous
space under the action of $\af \otimes_k \Ext^1(df_0, f_0)$;
\item if a lifting exists, then the automorphism group of liftings
of an element is canonically isomorphic to $\af \otimes_k \Ext^0(df_0, f_0)$.
\end{enumerate}
\end{theorem}

\proof cf. \cite{ZivRan}. \qed

In the following, we will also need results about deformation of
curves endowed with an invertible sheaf and a global section
of it. Before stating the result, we need to introduce some notation.

For any scheme $X/S$ denote by $\Pp^n_{X/S}$ the sheaf of 
principal part of order $n$ of $\O_X$ and for any invertible
sheaf on $X$, write $\Pp^n_{X/S}(\L)$ for the sheaf of principal
part of order $n$ of $\L$ (in particular, we have 
$\Pp^n_{X/S}(\L)=\Pp^n_{X/S} \otimes_{\O_X} \L$, the module structure on 
$\Pp^n_{X/S}$ being given by $d^n:\O_X \to \Pp^n_{X/S}$ 
(cf. \cite{EGA} IV, chapitre 16). In particular, one has
$\Ext^0_{\O_X}(\Pp^1_{X/S}(\L), \L)=\mathscr{D}iff(\L, \L)$.

\begin{theorem}
Let $k$ be a field, $X_0/k$ be a local complete intersection morphism
and $\L_0$ be an invertible sheaf on $X$. Let $A$ and $A'$ be local
artinian ring with residue field $k$, $A' \to A$ a small extension 
with kernel $\af$ and $(X, \L)$ be a deformation of $(X_0, \L_0)$ over 
$A$. Then the following is true
\begin{enumerate}
\item there exists a canonical element $\omega \in \af \otimes_k \Ext^2(\Pp^1_{X_0/k}(\L), \L)$ which vanishes if and only if there exists a lifting
of $(X, \L)$ to $A'$;
\item if there exists a lifting of $(X, \L)$ to $A'$, then the
set of isomorphism classes of liftings is naturally a principal 
homogeneous space under the action of $\af \otimes_k \Ext^1(\Pp^1_{X_0/k}(\L), \L)$;
\item if there exists a lifting $(X', \L')$ of $(X, \L)$ to $A'$ then 
the automorphism group of liftings of $(X', \L')$ is 
isomorphic to $\af \otimes_k \Ext^0(\Pp^1_{X/S}(\L), \L)$.
\end{enumerate}
\end{theorem}

\proof This is a generalization to the case of local complete
intersection morphism of the results of \cite{KodairaSpencer}
(see also \cite{GFGA}, page 13, corollaire 2). \qed

We can generalize results about the deformation of
schemes endowed with an invertible sheaf and a global section
of this sheaf.

\begin{proposition}
\label{AutLiftingSection}
Let $k$ be a field, $X_0/k$ a local complete intersection morphism,
$\L_0$ an invertible sheaf on $X_0$ and $s_0 \in \H^0(X_0, \L_0)$.
Then for any small extension $A'\to A$ with kernel $\af$ between local 
artinian rings with residue field $k$, any deformation $(X, \L, s)$
of $(X_0, \L_0, s_0)$ over $A$ and any lifting $(X', \L',s')$ of $(X, \L, s)$
to $A'$, the automorphism group of lifting of $(X', \L', s')$ is 
canonically isomorphic to 
$$\af \otimes_k \ker \left(\Ext^0_{\O_X}(\Pp^1_{X_0/k}(\L), \L) \overset{d^1_{s_0}}{\to} \H^0(X_0, \L_0))\right)$$
where $d^1_{s_0}(D)=Ds$.
\end{proposition}

\proof Follow the idea of Welters in \cite{Welters}, chapter 1. \qed

Next we need some results about the deformation of purely inseparable
morphism. 

\begin{lemma}
\label{StructureFrobenius}
Let $A$ be a complete local ring and $f:\Spec B\to \Spec C$ a finite
morphism between smooth $A$-curves which is purely inseparable of 
degree $p$ and suppose that $\Omega_{C/A}$ admits a basis which is 
exact. Then $f$ is of the form $C\to C[Y]/(P)$ the reduction of $P$
modulo the maximal ideal being of the form $Y^p-x$.

In particular, $\Omega_{X/A}= \bigoplus_{i=0}^{p-1} y^i C dy$.
\end{lemma}

\proof The general case is deduced from the case of a field using
flatness of $B$ over $A$. Hence we can assume $A$ to be a field 
so that $f$ can be identified to the Frobenius. Then choose a basis
$dx$ of $\Omega_{Y/A}$ which is exact. As $f$ is the 
Frobenius, $x$ is a $p$-th power in $B$, denote by $y \in B$
an element such that $y^p=x$. Since $dx$ is a basis of $\Omega_{C/k}$,
$\Spec C[y]/(y^p-x)$ is smooth over $k$, hence normal. Thus the 
result. \qed

\begin{proposition}
\label{LiftingFrobUpToDiff}
Let $A' \to A$ be a small extension of local artinian rings with same
residue field $k$, $f:X \to Y$ a morphism of smooth curves satisfying
the hypothesis of lemma \ref{StructureFrobenius} (in particular,
there exists a basis $dx$ of $\Omega_{Y/A}$). Let $X'$ and $Y'$ be 
liftings of $X$ and $Y$ respectively and $f'_1$ and $f'_2$ be two 
liftings of $f$ between $X'$ and $Y'$ such that
there exists $dx' \in \Omega_{Y'/A'}$ satisfying $df'_1(dx')=df'_2(dx')$.
Then $f'_1$ and $f'_2$ are equal up to an isomorphism of $Y'$.
\end{proposition}

\proof Let $\varepsilon$ a generator of the kernel $A' \to A$. Denote
by $\m'$ the maximal ideal of $A'$, $k=A'/\m'$ and $f_0:X_0 \to Y_0$
the special fiber of $f$.
As $f'_1$ and $f'_2$ are liftings of $f$, we have
$f'_1(x')-f'_2(x')=\varepsilon \mu$ with $\mu \in \O_{X_0}$.
By hypothesis, we have $\varepsilon d\mu=d(\varepsilon \mu)=0$.
Writing $\O_{X_0}=\bigoplus_{i=0}^{p-1} y^i \O_{Y_0}$ (such a decomposition
exists by lemma \ref{StructureFrobenius}) we see
that we have naturally $\mu \in \O_{Y_0}$.
Then we see that the automorphism of $Y'$ defined by
$$\varepsilon \mu \frac{\partial}{\partial x'} \in (\varepsilon) \otimes_k \Ext^0_{\O_{Y_0}}(\Omega_{Y_0/k}, \O_{Y_0})$$
sends $f'_1(x')$ to $f'_2(x')$, and thus sends $f'_1$ to $f'_2$.

\providecommand{\bysame}{\leavevmode\hbox to3em{\hrulefill}\thinspace}


\begin{thebibliography}{LMB00}

\bibitem[Bou83]{BourbakiAC2}
N.~Bourbaki, \emph{El\'ements de math\'ematique. alg\`ebre commutative. }, Masson, 1983.

\bibitem[DM69]{Deligne_Mumford}
P.~Deligne and D.~Mumford, \emph{The irreducibility of the space of curves of
  given genus}, Inst. Hautes \'Etudes Sci. Publ. Math.:36 (1969), 75--109.

\bibitem[EGA]{EGA}
A.~Grothendieck, \emph{\'{E}l\'ements de g\'eom\'etrie alg\'ebrique}, Inst.
  Hautes \'Etudes Sci. Publ. Math. \textbf{4, 8, 11, 17, 20, 24, 28, 32}
  (1961-1967).

\bibitem[Gro95]{GFGA}
Alexander Grothendieck, \emph{G\'eom\'etrie formelle et g\'eom\'etrie
  alg\'ebrique}, S\'eminaire Bourbaki, Vol.\ 5, Soc. Math. France, Paris, 1995,
  pp.~Exp.\ No.\ 182, 193--220, errata p.\ 390.

\bibitem[Hen00]{Henrio_p_adique}
Y.~Henrio, \emph{Arbres de hurwitz et automorphismes d'ordre $p$ des disques et
  des couronnes $p$-adiques formelles}, 2000, \mbox{arXiv:math.AG/0011098},
  preprint.

\bibitem[Hen01]{Henrio_galoisien}
Y.~Henrio, \emph{Rel\`evement galoisien des rev\^etements de courbes nodales},
  Manuscripta Math. \textbf{106}:2 (2001), 131--150.

\bibitem[HM82]{Harris_Mumford}
J.~Harris and D.~Mumford, \emph{On the {K}odaira dimension of the moduli space
  of curves}, Invent. Math. \textbf{67}:1 (1982), 23--88, With an appendix by
  William Fulton.

\bibitem[HS99]{HarbaterStevenson}
D.~Harbater and K.~F. Stevenson, \emph{Patching and thickening problems}, J.
  Algebra \textbf{212}:1 (1999), 272--304.

\bibitem[KS58]{KodairaSpencer}
K.~Kodaira and D.~C. Spencer, \emph{On deformations of complex analytic
  structures. {I}, {II}}, Ann. of Math. (2) \textbf{67} (1958), 328--466.

\bibitem[Liu]{LiuStableReduction}
Q.~Liu, \emph{Stable reduction of finite covers of curves}, preprint.

\bibitem[Liu02]{LiuBook}
Q.~Liu, \emph{Algebraic geometry and arithmetic curves}, Oxford University
  Press, 2002, Oxford Graduate Texts in Mathematics, No. 6.

\bibitem[LMB00]{Champs}
G.~Laumon and L.~Moret-Bailly, \emph{Champs alg\'ebriques}, Springer-Verlag,
  Berlin, 2000.

\bibitem[Mau03]{Papier1}
S.~Maugeais, \emph{Rel\`evement des rev\^etements p-cycliques des courbes
  rationnelles semi-stables}, Math. Ann. \textbf{327}:2 (2003), 365--393.

\bibitem[Moc95]{Mochizuki}
S.~Mochizuki, \emph{The geometry of the compactification of the {H}urwitz
  scheme}, Publ. Res. Inst. Math. Sci. \textbf{31}:3 (1995), 355--441.

\bibitem[Ran89]{ZivRan}
Ziv Ran, \emph{Deformations of maps}, Algebraic curves and projective geometry
  (Trento, 1988), Lecture Notes in Math., vol. 1389, Springer, 1989,
  pp.~246--253.

\bibitem[Sch68]{Deformation}
M.~Schlessinger, \emph{Functors of {A}rtin rings}, Trans. Amer. Math. Soc.
  \textbf{130} (1968), 208--222.

\bibitem[Wel83]{Welters}
Gerald~E. Welters, \emph{Polarized abelian varieties and the heat equations},
  Compositio Math. \textbf{49}:2 (1983), 173--194.

\bibitem[Wew98]{TheseWewers}
S~Wewers, \emph{Construction of hurwitz spaces}, Th\`ese de doctorat,
  Universität Essen, 1998.

\end{thebibliography}
\end{document}